\documentclass[11pt]{article}%
\usepackage[utf8]{inputenc}
\usepackage[T1]{fontenc}
\usepackage[english]{babel}
\usepackage{indentfirst}
\usepackage{latexsym}
\usepackage[colorlinks=true, bookmarksnumbered=true, bookmarksopen=true,
bookmarksopenlevel=3, pdfstartview=FitH, linkcolor=cyan, pdfmenubar=true,
pdftoolbar=true, bookmarks=true,citecolor=cyan, urlcolor=magenta,
filecolor=cyan,menucolor=black,plainpages=false,pdfpagelabels]{hyperref}
\usepackage[lmargin=2.4cm,tmargin=2.4cm,bmargin=2.4cm,rmargin=2.4cm]%
{geometry}
\usepackage{amsbsy, amsmath,amsfonts,amssymb,amsthm}
\usepackage{times}
\usepackage{graphicx}
\usepackage{amsmath}
\usepackage{amsfonts}
\usepackage{amssymb}%
\def\fin { \vskip 0pt \hfill $\diamond$ \vskip 12pt}
\newtheorem{theorem}{Theorem}[section]
\newtheorem{definition}[theorem]{Definition}

\newtheorem{lemma}[theorem]{Lemma}
\newtheorem{remark}[theorem]{Remark}

\numberwithin{equation}{section}
\begin{document}
\title{Besov-weak-Herz spaces and global solutions for Navier-Stokes equations}
\author{{{Lucas C. F. Ferreira}{\thanks{L. Ferreira was supported by FAPESP and CNPQ,
Brazil. Email:lcff@ime.unicamp.br (corresponding author).}} \ \ \ and \ \ {Jhean E. P\'{e}rez-L\'{o}pez \vspace{0.1cm}} {\thanks{J. P\'{e}rez-L\'{o}pez was supported by CAPES,
Brazil. Email:jhean754@ime.unicamp.br.}}}\\{\small Universidade Estadual de Campinas, IMECC - Departamento de
Matem\'{a}tica,} \\{\small {Rua S\'{e}rgio Buarque de Holanda, 651, CEP 13083-859, Campinas-SP,
Brazil.}}\vspace{-0.6cm}}
\date{}
\maketitle

\begin{abstract}
We consider the incompressible Navier-stokes equations (NS) in $\mathbb{R}%
^{n}$ for $n\geq2$. Global well-posedness is proved in critical
Besov-weak-Herz spaces (BWH-spaces) that consist in Besov spaces based on weak-Herz spaces.
These spaces are larger than some critical spaces considered in previous works
for (NS). For our purposes, we need to develop a basic theory for BWH-spaces containing properties and
estimates such as heat semigroup estimates, embedding
theorems, interpolation properties, among others. In particular, it is proved
a characterization of Besov-weak-Herz spaces as interpolation of
Sobolev-weak-Herz ones, which is key in our arguments. Self-similarity and
asymptotic behavior of solutions are also discussed. Our class of spaces and its properties developed here could also be employed to study other PDEs of elliptic, parabolic and conservation-law type.

\bigskip{} \noindent\textbf{AMS MSC:} 76D05; 76D03; 35A23; 35K08; 42B35;
46B70; 35C06; 35C15

\medskip{} \noindent\textbf{Key:} Navier-Stokes equations; Well-posedness;
Besov-weak-Herz spaces; Interpolation; Heat semigroup estimates; Self-similarity
\end{abstract}

\section{Introduction}

This paper is concerned with the incompressible Navier-Stokes equations%

\begin{equation}%
\begin{cases}
\frac{\partial u}{\partial t}-\Delta u+u\cdot\nabla u+\nabla\rho=0 & \text{in
}\,\,\mathbb{R}^{n}\times(0,\infty)\\
\nabla\cdot u=0 & \text{in}\,\,\mathbb{R}^{n}\times(0,\infty)\\
u(0)=u_{0} & \text{in}\,\,\mathbb{R}^{n}%
\end{cases}
, \label{eq:Problema NS introduccion}%
\end{equation}
where $n\geq2,$ $\rho$ is the pressure, $u=(u_{j})_{j=1}^{n}$ is the velocity
field and $u_{0}$ is a given initial velocity satisfying $\nabla\cdot u_{0}=0$.

After applying the Leray-Hopf projector $\mathbb{P}$ and using Duhamel's
principle, the Cauchy problem \eqref{eq:Problema NS introduccion} can be
reduced to the integral formulation
\begin{equation}
u(t)=G\left(  t\right)  u_{0}-\int_{0}^{t}G(t-\tau)\mathbb{P}\mathrm{div}%
\left(  u\otimes u\right)  (\tau)d\tau:=G\left(  t\right)  u_{0}+B(u,u)(t),
\label{eq:Integral NS}%
\end{equation}
where $u\otimes v:=\left(  u_{i}v_{j}\right)  _{1\leq i,j\leq n}$ is a
matrix-valued function and $G(t)=e^{t\Delta}$ is the heat semigroup. The
operator $\mathbb{P}$ can be expressed as $\mathbb{P}=(\mathbb{P}%
_{i,j})_{n\times n}$ where $\mathbb{P}_{i,j}:=\delta_{i,j}+\mathcal{R}%
_{i}\mathcal{R}_{j}$, $\delta_{i,j}$ is the Kronecker delta and $\mathcal{R}%
_{i}=(-\Delta)^{-1/2}\partial_{i}$ is the $i$-th Riesz transform.
Divergence-free solutions for \eqref{eq:Integral NS} are called mild solutions
for \eqref{eq:Problema NS introduccion}. Note that if $u$ is a smooth solution
for \eqref{eq:Problema NS introduccion} (or \eqref{eq:Integral NS}), then
\begin{equation}
u_{\lambda}(x,t):=\lambda u(\lambda x,\lambda^{2}t)
\label{eq:Scaling Navier-Stokes}%
\end{equation}
is also a solution with initial data
\begin{equation}
(u_{0})_{\lambda}(x)=\lambda u_{0}(\lambda x).
\label{eq:Scaling dado inicial NS}%
\end{equation}

Recall that given a Banach space $Y$ we say that it has scaling degree equal
to $k\in\mathbb{R}$ if $\left\Vert f(\lambda x)\right\Vert _{Y}\approx
\lambda^{k}\left\Vert f\right\Vert _{Y}$ for all $\lambda>0$ and $f\in Y$.
Motivated by \eqref{eq:Scaling dado inicial NS}, a Banach space $Y$ is called
critical for \eqref{eq:Problema NS introduccion} if it has scaling degree
equal to $-1$, that is, if $\left\Vert f\right\Vert _{Y}\approx\left\Vert
\lambda f(\lambda x)\right\Vert _{Y}$ for all $\lambda>0$ and $f\in Y$. In
turn, a solution of \eqref{eq:Problema NS introduccion} which is invariant by
the scaling \eqref{eq:Scaling Navier-Stokes}, i.e. $u=u_{\lambda},$ is called
a self-similar solution of \eqref{eq:Problema NS introduccion}. Note that in
order to obtain self-similar solutions, the initial data should be homogeneous
of degree $-1$.

Over the years, global-in-time well-posedness of small solutions for
\eqref{eq:Problema NS introduccion} in critical spaces has attracted the
interest of a number of authors. Without making a complete list, we mention
the works in homogeneous Sobolev $\dot{H}^{1/2}(\mathbb{R}^{3})$ \cite{FuKa},
Lebesgue $L^{n}(\mathbb{R}^{n})$ \cite{Ka1}, Marcinkiewicz $L^{n,\infty
}(\mathbb{R}^{n})$ \cite{Barraza,Yam}, Morrey $\mathcal{M}_{q}^{n}%
(\mathbb{R}^{n})$ \cite{Giga,Ka2,Taylor-1}, weak-Morrey $\mathcal{M}%
_{q,\infty}^{n}\left(  \mathbb{R}^{n}\right)  $
\cite{MiaoYuan,Lemarie-2,Ferreira-1}, $PM^{n-1}$-spaces \cite{CannKarch},
Besov $\dot{B}_{p,\infty}^{\frac{n}{p}-1}(\mathbb{R}^{n})$ for $p>n$
\cite{Cann1}, Fourier-Besov $F\dot{B}_{p,\infty}^{n-1-\frac{n}{p}}$
\cite{IwabTakada,KoniYone}, homogeneous weak-Herz spaces $W\dot{K}_{n,\infty
}^{0}\left(  \mathbb{R}^{n}\right)  $ \cite{Tsutsui-1}, Fourier-Herz
$\mathcal{B}_{r}^{-1}=F\dot{B}_{1,r}^{-1}$ with $r\in\lbrack1,2]$
\cite{CannWu,IwabTakada,LeiLin}, homogeneous Besov-Morrey $\mathcal{N}%
_{r,q,\infty}^{\frac{n}{r}-1}$ with $r>n$ \cite{KoYa,Maz}, and $BMO^{-1}$
\cite{KocTat}. The reader can find other examples in the nice review
\cite{Lemarie}. Up until now, we know that $BMO^{-1}$ and $\mathcal{N}%
_{r,1,\infty}^{\frac{n}{r}-1}$ are maximal critical spaces for
\eqref{eq:Problema NS introduccion} in the sense that it is not known a larger
critical space in which small solutions of \eqref{eq:Problema NS introduccion}
are globally well-posed.

The propose of this paper is to provide a new critical Besov type class for
global well-posedness of solutions for (\ref{eq:Problema NS introduccion}) by
assuming a smallness condition on initial data norms. Here we consider
homogeneous Besov-weak-Herz spaces $\dot{B}W\dot{K}_{p,q,r}^{\alpha,s}$ that
are a type of Besov space based on homogeneous weak-Herz spaces $W\dot
{K}_{p,q}^{\alpha}$. They are a natural extension of the spaces $BK_{p,q,r}%
^{\alpha,s}$ introduced in \cite{Xu} (see Definition
\ref{Def: Definicios espacios Besov-weak-Herz} in subsection 2.2). The Herz
space $K_{p,q}^{\alpha}$ was introduced by Herz in \cite{Herz-1} but his
definition is not appropriate for our purposes. Later, Johnson
\cite{Jonhson-1} obtained a characterization of the $K_{p,q}^{\alpha}$-norm in
terms of $L^{p}$-norms over annuli which is the base for the definition of the
spaces $W\dot{K}_{p,q}^{\alpha}$ in \cite{Tsutsui-1} and is the same one that
we use in the present paper. In order to achieve our aims, we need to develop
properties for $W\dot{K}_{p,q,r}^{\alpha}$ and $\dot{B}W\dot{K}_{p,q,r}%
^{\alpha,s}$-spaces such as H\"{o}lder inequality, estimates for convolution
operators, embedding theorems, interpolation properties, among others (see
Section 2). In particular, it is proved a characterization of Besov-weak-Herz
spaces in terms of interpolation of Sobolev-weak-Herz ones, which is key in
our arguments (see Lemma \ref{Besov-weak-Herz-Interp}). Moreover, we prove
estimates for the heat semigroup, as well as for the bilinear term $B(u,v)$ in
(\ref{eq:Integral NS}), in the context of $\dot{B}W\dot{K}_{p,q,r}^{\alpha,s}%
$-spaces. We also point out that these spaces and their basic theory developed here could be employed to study other PDEs of elliptic, parabolic and conservation-law type. It is worthy to observe that some arguments in this paper are
inspired by some of those in \cite{KoYa} that analyzed
(\ref{eq:Problema NS introduccion}) in Besov-Morrey spaces.

In what follows, we state our global well-posedness result.

\begin{theorem}
\label{Teo-main} Let $1\leq q\leq\infty,$ $\frac{n}{2}<p<\infty$ and
$0\leq\alpha<\min\{1-\frac{n}{2p},\frac{n}{2p}\}$. There exist $\epsilon>0$
and $\delta>0$ such that if $u_{0}\in\dot{B}W\dot{K}_{p,q,\infty}%
^{\alpha,\alpha+\frac{n}{p}-1}$ with $\nabla\cdot u_{0}=0$ and $\left\Vert
u_{0}\right\Vert _{\dot{B}W\dot{K}_{p,q,\infty}^{\alpha,\alpha+\frac{n}{p}-1}%
}\leq\delta$, then there exists a unique mild solution $u\in L^{\infty
}((0,\infty);\,\dot{B}W\dot{K}_{p,q,\infty}^{\alpha,\alpha+\frac{n}{p}-1})$
for (\ref{eq:Problema NS introduccion}) such that
\[
\left\Vert u\right\Vert _{X}:=\left\Vert u\right\Vert _{L^{\infty}%
((0,\infty);\dot{B}W\dot{K}_{p,q,\infty}^{\alpha,\alpha+\frac{n}{p}-1}%
)}+\mathop{\sup}\limits_{t>0}t^{\frac{1}{2}-\left(  \frac{\alpha}{2}+\frac
{n}{4p}\right)  }\left\Vert u\right\Vert _{W\dot{K}_{2p,2q}^{\alpha}}%
\leq2\epsilon.
\]
Moreover, $u(t)\overset{\ast}{\rightharpoonup}u_{0}$ in $\dot{B}%
_{\infty,\infty}^{-1}$, as $t\rightarrow0^{+}$, and solutions depend
continuously on initial data.
\end{theorem}

We have the continuous inclusions $L^{n}\subset L^{n,\infty}\subset W\dot{K}_{n,\infty}^{0}%
\subset\dot{B}W\dot{K}_{n,\infty,\infty}^{0,0}$ (see Lemmas
\ref{Lem: relaciones entre espacios y WHerz} and
\ref{Lem: Inclusiones Besov-herz con Herz}) and
\[
\dot{H}^{\frac{n}{2}-1}\subset L^{n}\subset\dot{B}_{p,\infty}^{\frac{n}{p}%
-1}\subset\dot{B}W\dot{K}_{p,\infty,\infty}^{0,\frac{n}{p}-1},\text{ for
}p\geq n\text{ (see Remark \ref{rem-inclusion-1}).}%
\]
So our initial data class extends those of some previous works; for instance,
the ones in \cite{FuKa,Ka1,Barraza,Cann1,Yam,Tsutsui-1}.

Notice that the parameter $s$ corresponds to the regularity index of the Besov
type space $\dot{B}W\dot{K}_{p,q,r}^{\alpha,s}$. Considering the family $\{\dot{B}%
W\dot{K}_{p,\infty,\infty}^{0,\frac{n}{p}-1}\}_{p>n/2}$, in the positive regularity range $n/2<p<n$ we
are dealing with spaces smaller than those with $p>n$, because of the
Sobolev embedding $\dot{B}W\dot{K}_{p_{2},\infty,\infty}^{0,\frac{n}{p_{2}}%
-1}\subset\dot{B}W\dot{K}_{p_{1},\infty,\infty}^{0,\frac{n}{p_{1}}-1}$ when
$p_{2}<p_{1}$ (see Lemma \ref{Sobolev-Lemma-1}). For $p=2$, $n=3$ and
$s=1/2$ (positive regularity), one can show $\dot{B}W\dot{K}%
_{2,\infty,\infty}^{0,1/2}\subset BMO^{-1}$ by using duality and
$\mathbb{A}_{\infty}$-atom decomposition. However, for $p>n$ (negative regularity and larger spaces) it is not clear for us whether there are inclusion relations between $\dot{B}W\dot{K}_{p,\infty,\infty}^{0,\frac{n}{p}-1}$ and $BMO^{-1}$
or between $\dot{B}W\dot{K}_{p,\infty,\infty}^{0,\frac{n}{p}-1}$ and
$\mathcal{N}_{r,1,\infty}^{\frac{n}{r}-1}$ with $r>n$. In this sense, our
result seems to give a new critical initial data class for existence of small
global mild solutions for \eqref{eq:Problema NS introduccion}. In any case, it
would be suitable to recall that well-posedness involves more properties than
only existence of solutions, namely existence, uniqueness, persistence, and continuous
dependence on initial data, which together characterize a
good behavior of the Navier-Stokes flow in the considered space.

We finish with some comments about self-similarity and asymptotic behavior of
solutions. It is not difficult to see that for $n\leq p<\infty$ the function
$f(x)=\left\vert x\right\vert ^{-1}$ belongs to $\dot{B}W\dot{K}%
_{p,\infty,\infty}^{0,\frac{n}{p}-1}$. So, the homogeneous Besov-weak-Herz
spaces (at least some of them) contain homogeneous functions of degree $-1$.
Thus, if one assumes further that the initial data $u_{0}$ is a homogeneous
vector field of degree $-1,$ then a standard procedure involving a Picard type
sequence gives that the solution obtained in Theorem \ref{Teo-main} is in fact
self-similar. Moreover, following some estimates and arguments in the proof of
Theorem \ref{Teo-main}, with an extra effort, it is possible to prove that if
we have $u_{0}$ and $v_{0}$ satisfying $\mathop{\lim}\limits_{t\rightarrow
\infty}\left\Vert G(t)\left(  u_{0}-v_{0}\right)  \right\Vert _{\dot{B}%
W\dot{K}_{p,q,\infty}^{\alpha,\alpha+\frac{n}{p}-1}}=0,$ then
\[
\mathop{\lim}\limits_{t\rightarrow\infty}\left\Vert u(\cdot,t)-v(\cdot
,t)\right\Vert _{\dot{B}W\dot{K}_{p,q,\infty}^{\alpha,\alpha+\frac{n}{p}-1}%
}=0,
\]
where $u$ and $v$ are the solutions obtained in Theorem \ref{Teo-main} with
initial data $u_{0}$ and $v_{0}$, respectively.

The plan of this paper is as follows. Section 2 is devoted to function spaces
where Herz and Sobolev-Herz spaces are considered in subsection 2.1 while
Sobolev-weak-Herz and Besov-weak-Herz spaces are addressed in subsection 2.2.
The proof of Theorem \ref{Teo-main} is performed in the final section through
three subsections, namely 3.1, 3.2 and 3.3. In the first we provide linear
estimates for the heat semigroup. The second is devoted to bilinear estimates
for $B(\cdot,\cdot)$ in our setting. After obtaining the needed estimates, the
proof is concluded in subsection 3.3 by means of a contraction argument.

\section{Function spaces}

In this section we recall some definitions and properties about function
spaces that will be considered throughout this paper.

\subsection{Weak-Herz and Sobolev weak-Herz spaces}

For an integer $k\in\mathbb{Z},$ we define the set $A_{k}$ as%

\begin{equation}
A_{k}=\left\{  x\in\mathbb{R}^{n};\,2^{k-1}\leq\left\vert x\right\vert
<2^{k}\right\}  , \label{decomp-1}%
\end{equation}
and observe that $\mathbb{R}^{n}\backslash\left\{  0\right\}  =\cup
_{k\in\mathbb{Z}}A_{k}.$ Taking $x\in A_{k}$ we have that%

\begin{align*}
y  &  \in A_{m}\,\mbox{and}\,\,m\leq k\,\Rightarrow2^{k-1}-2^{m}\leq\left\vert
x-y\right\vert <2^{k}+2^{m},\\
y  &  \in A_{m}\,\mbox{and}\,\,m\geq k\,\Rightarrow2^{m-1}-2^{k}\leq\left\vert
x-y\right\vert <2^{m}+2^{k}.
\end{align*}
Consider also the sets%

\begin{align}
C_{m,k}  &  =\left\{  \xi;2^{k-1}-2^{m}\leq\left\vert \xi\right\vert
<2^{k}+2^{m}\right\}  ,\nonumber\\
\widetilde{C}_{m,k}  &  =\left\{  \xi;2^{m-1}-2^{k}\leq\left\vert
\xi\right\vert <2^{m}+2^{k}\right\}  . \label{eq:Cmklinha}%
\end{align}

Now we are able to define the weak-Herz spaces.

\begin{definition}
\label{Definition Weak-herz spaces} Let $1<p\leq\infty$, $1\leq q\leq\infty$
and $\alpha\in\mathbb{R}$. The Homogeneous weak-Herz space $W\dot{K}%
_{p,q}^{\alpha}=$ $W\dot{K}_{p,q}^{\alpha}(\mathbb{R}^{n})$ is defined as the
set of all measurable functions such that the following quantity is finite
\begin{equation}
\left\Vert f\right\Vert _{W\dot{K}_{p,q}^{\alpha}}:=\left\{
\begin{array}
[c]{l}%
\left(  \sum\limits_{k\in\mathbb{Z}}2^{k\alpha q}\left\Vert f\right\Vert
_{L^{p,\infty}(A_{k})}^{q}\right)  ^{\frac{1}{q}}\,\ \text{if }q<\infty,\\
\mathop{\sup}\limits_{k\in\mathbb{Z}}\,2^{k\alpha}\left\Vert f\right\Vert
_{L^{p,\infty}(A_{k})}\,\,\,\,\,\,\,\,\,\,\,\ \ \,\text{if \ }q=\infty.
\end{array}
\right.  \label{norm-Herz-1}%
\end{equation}

\end{definition}

For $\alpha\in\mathbb{R},$ $1<p\leq\infty$ and $1\leq q\leq\infty$, the
quantity $\left\Vert \cdot\right\Vert _{W\dot{K}_{p,q}^{\alpha}}$ defines a
norm in $W\dot{K}_{p,q}^{\alpha}$ and the pair $(W\dot{K}_{p,q}^{\alpha
},\left\Vert \cdot\right\Vert _{W\dot{K}_{p,q}^{\alpha}})$ is a Banach space
(see e.g. \cite{Hernandez-Yang,Tsutsui-1}).

H\"{o}lder inequality holds in the setting of homogeneous Weak-Herz spaces
(see \cite{Tsutsui-1}). To be more precise, if $1<p,p_{1},p_{2}\leq\infty,$
$1\leq q,q_{1},q_{2}\leq\infty$ and $\alpha,\alpha_{1},\alpha_{2}\in
\mathbb{R}$ are such that $\frac{1}{p}=\frac{1}{p_{1}}+\frac{1}{p_{2}},$
$\frac{1}{q}=\frac{1}{q_{1}}+\frac{1}{q_{2}},$ and $\alpha=\alpha_{1}%
+\alpha_{2}$, then%

\begin{equation}
\left\Vert fg\right\Vert _{W\dot{K}_{p,q}^{\alpha}}\leq C\left\Vert
f\right\Vert _{W\dot{K}_{p_{1},q_{1}}^{\alpha_{1}}}\left\Vert g\right\Vert
_{W\dot{K}_{p_{2},q_{2}}^{\alpha_{2}}}, \label{eq:Desigualdade de Holder.}%
\end{equation}
where $C>0$ is an universal constant. In fact, for all $k\in\mathbb{Z},$ we have%

\[
\left\Vert fg\right\Vert _{L^{p,\infty}(A_{k})}\leq C\left\Vert f\right\Vert
_{L^{p_{1},\infty}(A_{k})}\left\Vert g\right\Vert _{L^{p_{2},\infty}(A_{k})},
\]
and therefore%

\begin{align}
\left\Vert fg\right\Vert _{W\dot{K}_{p,q}^{\alpha}}  &  =\left(
\sum\limits_{k\in\mathbb{Z}}2^{k\alpha q}\left\Vert fg\right\Vert
_{L^{p,\infty}(A_{k})}^{q}\right)  ^{1/q}\nonumber\\
&  \leq C\left(  \sum\limits_{k\in\mathbb{Z}}2^{k\alpha_{1}q}\left\Vert
f\right\Vert _{L^{p_{1},\infty}(A_{k})}^{q}2^{k\alpha_{2}q}\left\Vert
g\right\Vert _{L^{p_{2},\infty}(A_{k})}^{q}\right)  ^{1/q}\nonumber\\
&  \leq C\left\Vert f\right\Vert _{W\dot{K}_{p_{1},q_{1}}^{\alpha_{1}}%
}\left\Vert g\right\Vert _{W\dot{K}_{p_{2},q_{2}}^{\alpha_{2}}}.
\label{aux-holder-1}%
\end{align}
Taking in particular $(\alpha_{1},p_{1},q_{1})=(0,\infty,\infty)$ in
(\ref{aux-holder-1}), we obtain%

\begin{equation}
\left\Vert fg\right\Vert _{W\dot{K}_{p,q}^{\alpha}}\leq C\left\Vert
f\right\Vert _{L^{\infty}(\mathbb{R}^{n})}\left\Vert g\right\Vert _{W\dot
{K}_{p,q}^{\alpha}}. \label{eq:Holder con Linfinito}%
\end{equation}

More below we will need to estimate some convolution operators, particularly
the heat semigroup, in weak-Herz and Besov-weak-Herz spaces. The following
lemma will be useful for that propose.

\begin{lemma}
\label{lem:[Comvolucao]} (Convolution) Let $1\leq p_{1}<\infty$ and
$1<r,p_{2}<\infty$ be such that $1+\frac{1}{r}=\frac{1}{p_{1}}+\frac{1}{p_{2}%
}$. Also, let $1\leq q\leq\infty,$ $-\frac{n}{r}<\alpha<n\left(  1-\frac
{1}{p_{2}}\right)  $, and $\theta\in L^{p_{1}}\left(  \mathbb{R}^{n}\right)  $
be such that $\theta\left\vert \cdot\right\vert ^{n/p_{1}}\in L^{\infty
}\left(  \mathbb{R}^{n}\right)  .$ Then, there exists a positive constant $C$
independent of $\theta$ such that
\begin{equation}
\left\Vert \theta\ast f\right\Vert _{W\dot{K}_{r,q}^{\alpha}}\leq
C\max\left\{  \left\Vert \theta\right\Vert _{L^{p_{1}}},\left\Vert \left\vert
\cdot\right\vert ^{n/p_{1}}\theta\right\Vert _{L^{\infty}}\right\}  \left\Vert
f\right\Vert _{W\dot{K}_{p_{2},q}^{\alpha}}, \label{aux-conv-1-0}%
\end{equation}
for all $f\in W\dot{K}_{p_{2},q}^{\alpha}$.
\end{lemma}

\textbf{Proof. }Denote $f_{m}=f|_{A_{m}}.$ Recalling the decomposition
(\ref{decomp-1}), for $k\in\mathbb{Z}$ we can estimate
\begin{align}
&  2^{k\alpha}\left\Vert \theta\ast f\right\Vert _{L^{r,\infty}\left(
A_{k}\right)  }\nonumber\\
&  \leq2^{k\alpha}\left\{  \left\Vert \sum\limits_{m\leq k-2}\theta\ast
f_{m}\right\Vert _{L^{r,\infty}\left(  A_{k}\right)  }+\left\Vert
\sum\limits_{m=k-1}^{k+1}\theta\ast f_{m}\right\Vert _{L^{r,\infty}\left(
A_{k}\right)  }+\left\Vert \sum\limits_{m\geq k+2}\theta\ast f_{m}\right\Vert
_{L^{r,\infty}\left(  A_{k}\right)  }\right\} \nonumber\\
&  :=I_{1}^{k}+I_{2}^{k}+I_{3}^{k}. \label{aux-conv-1}%
\end{align}

Using the notations in (\ref{eq:Cmklinha}) and the change of variable $z=k-m$,
we handle the term $I_{3}^{k}$ as follows%

\begin{align}
I_{3}^{k}  &  \leq2^{k\alpha}\left\Vert \sum\limits_{m\geq k+2}\theta\ast
f_{m}\right\Vert _{L^{r,\infty}\left(  A_{k}\right)  }\leq2^{k\alpha
}\left\Vert \sum\limits_{m\geq k+2}\theta\ast f_{m}\right\Vert _{L^{r}\left(
A_{k}\right)  }\nonumber\\
&  \leq2^{k\alpha}\left(  \int_{A_{k}}\left\vert \sum\limits_{m\geq k+2}%
\int_{\mathbb{R}^{n}}\theta(x-y)f_{m}(y)dy\right\vert ^{r}dx\right)
^{1/r}\nonumber\\
&  =2^{k\alpha}\left(  \int_{A_{k}}\left\vert \sum\limits_{m\geq k+2}%
\int_{\mathbb{R}^{n}}\theta(x-y)\chi_{\tilde{C}_{m.k}}(x-y)f_{m}%
(y)dy\right\vert ^{r}dx\right)  ^{1/r}\nonumber\\
&  \leq C\left\Vert \left\vert \cdot\right\vert ^{n/p_{1}}\theta\right\Vert
_{L^{\infty}}2^{k\alpha}\left(  \int_{A_{k}}\left(  \sum\limits_{m\geq
k+2}\int_{\mathbb{R}^{n}}\left\vert x-y\right\vert ^{-n/p_{1}}\chi_{\tilde
{C}_{m.k}}(x-y)\left\vert f_{m}(y)\right\vert dy\right)  ^{r}dx\right)
^{1/r}\nonumber\\
&  \leq C\left\Vert \left\vert \cdot\right\vert ^{n/p_{1}}\theta\right\Vert
_{L^{\infty}}2^{k\alpha}\left(  \int_{A_{k}}\left(  \sum\limits_{m\geq
k+2}2^{-mn/p_{1}}\left\Vert f\right\Vert _{L^{1}\left(  A_{m}\right)
}\right)  ^{r}dx\right)  ^{1/r}. \label{aux-conv-est-001}%
\end{align}
Recalling the inclusion $L^{p_{2},\infty}\left(  A_{m}\right)  \hookrightarrow
$ $L^{1}\left(  A_{m}\right)  ,$ we can continue to estimate the right-hand
side of (\ref{aux-conv-est-001}) in order to obtain%

\begin{align}
&  \text{R.H.S. of (\ref{aux-conv-est-001})}\nonumber\\
&  \leq C\left\Vert \left\vert \cdot\right\vert ^{n/p_{1}}\theta\right\Vert
_{L^{\infty}}2^{k\alpha}\left(  \int_{A_{k}}\left(  \sum\limits_{m\geq
k+2}2^{-mn/p_{1}}2^{nm\left(  1-1/p_{2}\right)  }\left\Vert f\right\Vert
_{L^{p_{2},\infty}\left(  A_{m}\right)  }\right)  ^{r}dx\right)
^{1/r}\nonumber\\
&  \leq C\left\Vert \left\vert \cdot\right\vert ^{n/p_{1}}\theta\right\Vert
_{L^{\infty}}2^{k\alpha}2^{k\frac{n}{r}}\sum\limits_{m\geq k+2}2^{-m\frac
{n}{r}}\left\Vert f\right\Vert _{L^{p_{2},\infty}\left(  A_{m}\right)
}\nonumber\\
&  \leq C\left\Vert \left\vert \cdot\right\vert ^{n/p_{1}}\theta\right\Vert
_{L^{\infty}}\sum\limits_{-2\geq z}2^{k\left(  \alpha+\frac{n}{r}\right)
}2^{(z-k)\frac{n}{r}}\left\Vert f\right\Vert _{L^{p_{2},\infty}\left(
A_{k-z}\right)  }\nonumber\\
&  \leq C\left\Vert \left\vert \cdot\right\vert ^{n/p_{1}}\theta\right\Vert
_{L^{\infty}}\sum\limits_{-2\geq z}2^{k\alpha}2^{z\frac{n}{r}}2^{-\left(
k-z\right)  \alpha}2^{\left(  k-z\right)  \alpha}\left\Vert f\right\Vert
_{L^{p_{2},\infty}\left(  A_{k-z}\right)  }\nonumber\\
&  \leq C\left\Vert \left\vert \cdot\right\vert ^{n/p_{1}}\theta\right\Vert
_{L^{\infty}}\sum\limits_{-2\geq z}2^{z\left(  \frac{n}{r}+\alpha\right)
}2^{\left(  k-z\right)  \alpha}\left\Vert f\right\Vert _{L^{p_{2},\infty
}\left(  A_{k-z}\right)  }. \label{aux-conv-est-1}%
\end{align}
The above estimates and Minkowski inequality lead us to (with usual
modification in the case $q=\infty$)%

\[
\left(  \sum\limits_{k\in\mathbb{Z}}\left(  I_{3}^{k}\right)  ^{q}\right)
^{1/q}\leq CM_{\theta}\left\Vert f\right\Vert _{W\dot{K}_{p_{2},q}^{\alpha}}.
\]

For the parcel $I_{2}^{k}$, we estimate%

\begin{align*}
I_{2}^{k}  &  \leq2^{k\alpha}\sum\limits_{m=k-1}^{k+1}\left\Vert \theta\ast
f_{m}\right\Vert _{L^{r,\infty}\left(  A_{k}\right)  }\leq2^{k\alpha}%
\sum\limits_{m=k-1}^{k+1}\left\Vert \theta\ast f_{m}\right\Vert _{L^{r,\infty
}\left(  \mathbb{R}^{n}\right)  }\\
&  \leq2^{k\alpha}\sum\limits_{m=k-1}^{k+1}\left\Vert \theta\right\Vert
_{L^{p_{1},\infty}}\left\Vert f_{m}\right\Vert _{L^{p_{2},\infty}}\leq
C\left\Vert \theta\right\Vert _{L^{p_{1}}}\sum\limits_{l=-1}^{1}2^{\left(
k+l\right)  \alpha}\left\Vert f\right\Vert _{L^{p_{2},\infty}\left(
A_{k+l}\right)  },
\end{align*}
which implies%
\[
\left(  \sum\limits_{k\in\mathbb{Z}}\left(  I_{2}^{k}\right)  ^{q}\right)
^{1/q}\leq CM_{\theta}\left\Vert f\right\Vert _{W\dot{K}_{p_{2},q}^{\alpha}}.
\]

Proceeding similarly to the estimates (\ref{aux-conv-est-001}%
)-(\ref{aux-conv-est-1}) but considering $C_{m.k}$ in place of $\tilde
{C}_{m.k}$ , the parcel $I_{1}^{k}$ can be estimated as%

\begin{align}
I_{1}^{k}  &  \leq C\left\Vert \left\vert \cdot\right\vert ^{n/p_{1}}%
\theta\right\Vert _{L^{\infty}}2^{k\alpha}\left(  \int_{A_{k}}\left(
\sum\limits_{m\leq k-2}\int_{\mathbb{R}^{n}}\left\vert x-y\right\vert
^{-n/p_{1}}\chi_{C_{m.k}}(x-y)\left\vert f_{m}(y)\right\vert dy\right)
^{r}dx\right)  ^{1/r}\nonumber\\
&  \leq C\left\Vert \left\vert \cdot\right\vert ^{n/p_{1}}\theta\right\Vert
_{L^{\infty}}2^{k\alpha}\left(  \int_{A_{k}}\left(  \sum\limits_{m\leq
k-2}2^{-kn/p_{1}}\left\Vert f\right\Vert _{L^{1}\left(  A_{m}\right)
}\right)  ^{r}dx\right)  ^{1/r}\nonumber\\
&  \leq C\left\Vert \left\vert \cdot\right\vert ^{n/p_{1}}\theta\right\Vert
_{L^{\infty}}2^{k\alpha}2^{k\frac{n}{r}}\sum\limits_{m\leq k-2}2^{-kn/p_{1}%
}\left\Vert f\right\Vert _{L^{1}\left(  A_{m}\right)  }\nonumber\\
&  \leq C\left\Vert \left\vert \cdot\right\vert ^{n/p_{1}}\theta\right\Vert
_{L^{\infty}}\sum\limits_{m\leq k-2}2^{k\left(  \alpha-n+\frac{n}{p_{2}%
}\right)  }\left\Vert f\right\Vert _{L^{1}\left(  A_{m}\right)  }\nonumber\\
&  \leq C\left\Vert \left\vert \cdot\right\vert ^{n/p_{1}}\theta\right\Vert
_{L^{\infty}}\sum\limits_{2\leq z}2^{k\left(  \alpha-n+\frac{n}{p_{2}}\right)
}2^{n\left(  k-z\right)  \left(  1-1/p_{2}\right)  }\left\Vert f\right\Vert
_{L^{p_{2},\infty}\left(  A_{k-z}\right)  }\nonumber\\
&  \leq C\left\Vert \left\vert \cdot\right\vert ^{n/p_{1}}\theta\right\Vert
_{L^{\infty}}\sum\limits_{2\leq z}2^{z\left(  \alpha-n+\frac{n}{p_{2}}\right)
}2^{\left(  k-z\right)  \alpha}\left\Vert f\right\Vert _{L^{p_{2},\infty
}\left(  A_{k-z}\right)  }. \label{aux-conv-est-2}%
\end{align}
It follows from (\ref{aux-conv-est-2}) that%

\[
\left(  \sum\limits_{k\in\mathbb{Z}}\left(  I_{1}^{k}\right)  ^{q}\right)
^{1/q}\leq CM_{\theta}\left\Vert f\right\Vert _{W\dot{K}_{p_{2},q}^{\alpha}}.
\]
Finally, the desired estimate is obtained after recalling the norm
(\ref{norm-Herz-1}) and using the above estimates for $I_{j}^{k}$ in
(\ref{aux-conv-1}).

\fin

Let $\varphi\in C_{c}^{\infty}\left(  \mathbb{R}^{n}\backslash\left\{
0\right\}  \right)  $ be radially symmetric and such that%

\[
\mbox{supp}\,\varphi\subset\left\{  x\,;\,\frac{3}{4}\leq\left\vert
x\right\vert \leq\frac{8}{3}\right\}
\]
and%

\[
\sum\limits_{j\in\mathbb{N}}\varphi_{j}(\xi)=1,\,\text{\ }\forall\xi
\in\mathbb{R}^{n}\backslash\left\{  0\right\}  ,
\]
where $\varphi_{j}(\xi):=\varphi\left(  \xi2^{-j}\right)  .$ Now we can define
the well-known localization operators $\Delta_{j}$ and $S_{j}$%

\begin{align*}
\Delta_{j}f  &  =\varphi_{j}(D)f=(\mathcal{F}^{-1}\varphi_{j})\ast f,\\
S_{k}f  &  =\sum\limits_{j\leq k}\Delta_{j}f.
\end{align*}
It is easy to see that we have the identities%

\[
\Delta_{j}\Delta_{k}f=0\,\,\mbox{if}\,\,\left\vert j-k\right\vert \geq2\text{
and }\Delta_{j}\left(  S_{k-2}g\Delta_{k}f\right)
=0\,\,\mbox{if}\,\,\left\vert j-k\right\vert \geq5.
\]
Finally, Bony's decomposition (see \cite{Bony}) gives%

\[
fg=T_{f}g+T_{g}f+R(fg),
\]
where%

\[
T_{f}g=\sum\limits_{j\in Z}S_{j-2}f\Delta_{j}g,\text{ }R(fg)=\sum\limits_{j\in
Z}\Delta_{j}f\tilde{\Delta}_{j}g\text{ and }\tilde{\Delta}_{j}g=\sum
\limits_{\left\vert j-j^{\prime}\right\vert \leq1}\Delta_{j^{\prime}}g.
\]

The next lemma will be useful in order to estimate some multiplier operators
in Besov-Weak-Herz spaces.

\begin{lemma}
\label{Teo:Operador-pseudo-Weak-Herz} Let $1<p<\infty,$ $1\leq q\leq\infty,$
$-\frac{n}{p}<\alpha<n\left(  1-\frac{1}{p}\right)  ,$ $m\in\mathbb{R}$ and
$D_{j}=\left\{  x\,;\,\frac{3}{4}2^{j}\leq\left\vert x\right\vert \leq\frac
{8}{3}2^{j}\right\}  $ for $j\in\mathbb{Z}$. Let $P$ be a $C^{n}$-function on
$\tilde{D}_{j}:=D_{j-1}\cup D_{j}\cup D_{j+1}$ such that $\left\vert
\partial_{\xi}^{\beta}P\left(  \xi\right)  \right\vert \leq C2^{\left(
m-\left\vert \beta\right\vert \right)  j}$ for all $\xi\in\tilde{D}_{j}$ and
multi-index $\beta$ satisfying $\left\vert \beta\right\vert \leq\left[
n/2\right]  +1$. Then, we have that%
\[
\left\Vert \left(  P\hat{f}\right)  ^{\vee}\right\Vert _{W\dot{K}%
_{p,q}^{\alpha}}\leq C2^{jm}\left\Vert f\right\Vert _{W\dot{K}_{p,q}^{\alpha}%
},
\]
for all $f\in W\dot{K}_{p,q}^{\alpha}$ such that $\mbox{supp}\hat{f}\subset
D_{j}$.
\end{lemma}

\textbf{Proof.} We start by defining $\tilde{\varphi}_{j}=\varphi
_{j-1}+\varphi_{j}+\varphi_{j+1}$ and $K(x)=\left(  P\tilde{\varphi}%
_{j}\right)  ^{\vee}$. Since $\mbox{supp}\hat{f}\subset D_{j}$ we have that
$P\left(  \xi\right)  \hat{f}\left(  \xi\right)  =P\left(  \xi\right)
\tilde{\varphi}_{j}(\xi)\hat{f}\left(  \xi\right)  $, and therefore $\left(
P\hat{f}\right)  ^{\vee}=\left(  P\tilde{\varphi}_{j}\hat{f}\right)  ^{\vee
}=K\ast f$.

Using Lemma \ref{lem:[Comvolucao]} we get%

\[
\left\Vert \left(  P\hat{f}\right)  ^{\vee}\right\Vert _{W\dot{K}%
_{p,q}^{\alpha}}\leq C\max\left\{  \left\Vert K\right\Vert _{L^{1}},\left\Vert
\left\vert \cdot\right\vert ^{n}K\right\Vert _{L^{\infty}}\right\}  \left\Vert
f\right\Vert _{W\dot{K}_{p,q}^{\alpha}}.
\]
It remains to show that $\max\left\{  \left\Vert K\right\Vert _{L^{1}%
},\left\Vert \left\vert \cdot\right\vert ^{n}K\right\Vert _{L^{\infty}%
}\right\}  \leq C2^{mj}$. For that, let $N\in\mathbb{N}$ be such that
$\frac{n}{2}<N\leq n$ and proceed as follows%

\begin{align*}
&  \left\Vert K\right\Vert _{L^{1}}=\int\limits_{B\left(  0,2^{-j}\right)
}K\left(  y\right)  +\int\limits_{\left\vert y\right\vert \geq2^{-j}}K\left(
y\right) \\
&  \leq\left(  \int_{B\left(  0,2^{-j}\right)  }1\right)  ^{1/2}\left(
\int_{B\left(  0,2^{-j}\right)  }\left\vert K\left(  y\right)  \right\vert
^{2}\right)  ^{1/2}+\left(  \int_{\left\vert y\right\vert \geq2^{-j}%
}\left\vert y\right\vert ^{-2N}\right)  ^{1/2}\left(  \int_{\left\vert
y\right\vert \geq2^{-j}}\left\vert y\right\vert ^{2N}\left\vert K\left(
y\right)  \right\vert ^{2}\right)  ^{1/2}\\
&  \leq C2^{-j\frac{n}{2}}\left\Vert P\tilde{\varphi}_{j}\right\Vert _{L^{2}%
}+C2^{-j\left(  -N+\frac{n}{2}\right)  }\sum\limits_{\left\vert \beta
\right\vert =N}\left\Vert \left(  \cdot\right)  ^{\beta}K\right\Vert _{L^{2}%
}\\
&  \leq C2^{-j\frac{n}{2}}\left\Vert P\tilde{\varphi}_{j}\right\Vert _{L^{2}%
}+C2^{-j\left(  -N+\frac{n}{2}\right)  }\sum\limits_{\left\vert \beta
\right\vert =N}\left\Vert \partial^{\beta}\left(  P\tilde{\varphi}_{j}\right)
\right\Vert _{L^{2}}\\
&  \leq C2^{-j\frac{n}{2}}C2^{mj}2^{j\frac{n}{2}}+C2^{-j\left(  -N+\frac{n}%
{2}\right)  }C2^{j\left(  m-N\right)  }2^{j\frac{n}{2}}\\
&  \leq C2^{mj}.
\end{align*}
For the norm $\left\Vert \left\vert \cdot\right\vert ^{n}K\right\Vert
_{L^{\infty}}$, we have that%

\begin{align*}
&  \left\Vert \left\vert \cdot\right\vert ^{n}K\right\Vert _{L^{\infty}}%
\leq\sum\limits_{\left\vert \beta\right\vert =n}\left\Vert \left(
\cdot\right)  ^{\beta}K\right\Vert _{L^{\infty}}\leq C\sum\limits_{\left\vert
\beta\right\vert =n}\left\Vert \partial^{\beta}\left(  P\tilde{\varphi}%
_{j}\right)  \right\Vert _{L^{1}}\\
&  \leq C\sum\limits_{\left\vert \beta\right\vert =n}2^{j\left(  m-n\right)
}2^{jn}\leq C2^{mj},
\end{align*}
as required.\fin

\subsection{Sobolev-weak-Herz spaces and Besov-weak-Herz spaces}

In this section we introduce the homogeneous Sobolev-weak-Herz spaces and
Besov-weak-Herz spaces. We also shall prove a number of properties about these
spaces that will be useful in our study of the Navier-Stokes equations. These
spaces are a generalization of Sobolev-Herz and Besov-Herz spaces found in
\cite{Xu}.

\begin{definition}
\label{Def: Definicios espacios Sobolev-weak-herz}Let $1<p\leq\infty,$ $1\leq
q\leq\infty$ and $\alpha,s\in\mathbb{R}$. Recall the Riesz operator
$\widehat{I^{s}f}=\left\vert \xi\right\vert ^{s}\hat{f}$. The homogeneous
Sobolev-weak-Herz spaces $W\dot{K}_{p,q}^{\alpha,s}=W\dot{K}_{p,q}^{\alpha
,s}(\mathbb{R}^{n})$ are defined as%
\begin{equation}
W\dot{K}_{p,q}^{\alpha,s}=\left\{  f\in\mathcal{S}^{\prime}(\mathbb{R}%
^{n})/\mathcal{P};\,\left\Vert I^{s}f\right\Vert _{W\dot{K}_{p,q}^{\alpha}%
}<\infty\right\}  . \label{eq:Definicion espacios de Herz com s}%
\end{equation}

\end{definition}

\begin{definition}
\label{Def: Definicios espacios Besov-weak-Herz} Let $1<p\leq\infty,$ $1\leq
q,r\leq\infty$ and $\alpha,s\in\mathbb{R}$. The homogeneous Besov-weak-Herz
spaces $\dot{B}W\dot{K}_{p,q,r}^{\alpha,s}=\dot{B}W\dot{K}_{p,q,r}^{\alpha
,s}\left(  \mathbb{R}^{n}\right)  $ are defined as%
\[
\dot{B}W\dot{K}_{p,q,r}^{\alpha,s}=\left\{  f\in\mathcal{S}^{\prime
}(\mathbb{R}^{n})/\mathcal{P};\,\left\Vert f\right\Vert _{\dot{B}W\dot
{K}_{p,q,r}^{\alpha,s}}<\infty\right\}  ,
\]
where%
\begin{equation}
\left\Vert f\right\Vert _{\dot{B}W\dot{K}_{p,q,r}^{\alpha,s}}:=\left\{
\begin{array}
[c]{l}%
\left(  \sum\limits_{j\in\mathbb{Z}}2^{jsr}\left\Vert \Delta_{j}f\right\Vert
_{W\dot{K}_{p,q}^{\alpha}}^{r}\right)  ^{1/r}\,\,\mbox{if}\,\,r<\infty\\
\mathop{\sup}\limits_{j\in\mathbb{Z}}\,2^{js}\left\Vert \Delta_{j}f\right\Vert
_{W\dot{K}_{p,q}^{\alpha}}%
\,\,\,\,\,\,\,\,\,\,\,\,\,\,\,\,\,\,\,\mbox{if}\,\,r=\infty.
\end{array}
\right.  \label{norm-1}%
\end{equation}

\end{definition}

\begin{remark}
\label{rem-inclusion-1}

\begin{itemize}
\item[(i)] The spaces $W\dot{K}_{p,q}^{\alpha,s}$ and $\dot{B}W\dot{K}%
_{p,q,r}^{\alpha,s}$ are Banach spaces endowed with the norms $\left\Vert
\cdot\right\Vert _{W\dot{K}_{p,q}^{\alpha}}$ and $\left\Vert \cdot\right\Vert
_{\dot{B}W\dot{K}_{p,q,r}^{\alpha,s}}$, respectively.

\item[(ii)] The continuous inclusion $\dot{B}_{p,r}^{s}\left(  \mathbb{R}%
^{n}\right)  \subset\dot{B}W\dot{K}_{p,\infty,r}^{0,s}\left(  \mathbb{R}%
^{n}\right)  $ holds for all $s\in\mathbb{R}$, $1<p\leq\infty$, and $1\leq
r\leq\infty,$ where $\dot{B}_{p,r}^{s}$ stands for homogeneous Besov spaces.
For that, it is sufficient to recall the definition of Besov spaces (see
\cite[pg.146]{BL}) and (\ref{norm-1}) and to use the inclusion $L^{p}\subset
W\dot{K}_{p,\infty}^{0}$ that is going to be showed in the lemma below.
\end{itemize}
\end{remark}

The next lemma contains relations between weak-$L^{p}$, weak-Herz and Morrey
spaces. For the definition and some properties about Morrey spaces we refer
the reader to \cite{KoYa} (see also \cite{Ka2} for an equivalent definition
and further properties).

\begin{lemma}
\label{Lem: relaciones entre espacios y WHerz} For $1<p<\infty$, we have the
continuous inclusion%
\begin{equation}
L^{p}\varsubsetneq L^{p,\infty}\varsubsetneq W\dot{K}_{p,\infty}^{0}.
\label{eq:Inclusao estrita Lp}%
\end{equation}

Moreover, let $\mathcal{M}_{q}^{r}$ stand for homogeneous Morrey spaces,
$1\leq q\leq r<\infty$ and $n/r\neq\alpha+n/p$ when $q<p$. Then
\begin{equation}
W\dot{K}_{p,\infty}^{\alpha}\not \subset \mathcal{M}_{q}^{r}.
\label{eq:No inclusion en Morrey}%
\end{equation}

\end{lemma}

\textbf{Proof.} The first inclusion in (\ref{eq:Inclusao estrita Lp}) is
well-known, so we only prove the second one. For that, it is sufficient to
note that $\left\Vert f\right\Vert _{L^{p,\infty}(A_{k})}\leq\left\Vert
f\right\Vert _{L^{p,\infty}(\mathbb{R}^{n})}$ $\forall k\in\mathbb{Z}$ and
after to take the supremum over $k$. In order to see the strictness of the
inclusion, take $x_{k}=\frac{3}{2}2^{k-1}\vec{e}_{1}$ and $h(x):=\sum
\limits_{k=1}^{\infty}\left\vert x-x_{k}\right\vert ^{-\frac{n}{p}}%
\chi_{B(0,1/8)}\left(  x-x_{k}\right)  .$ It is clear that $h\in W\dot
{K}_{p,\infty}^{0}$ but not to $L^{p,\infty}\left(  \mathbb{R}^{n}\right)  $.

Now we turn to (\ref{eq:No inclusion en Morrey}). For $f(x)=\left\vert
x\right\vert ^{-\frac{n}{p}}$, we have that $f\in L^{p,\infty}\subset W\dot
{K}_{p,\infty}^{0}$. On the other hand, for any $q\geq p$ note that
$\left\Vert f\right\Vert _{L^{q}\left(  B\left(  0,R\right)  \right)  }%
=\infty$, and then $f\notin\mathcal{M}_{q}^{r}$ for any $r$. Finally, if
$n/r\neq\alpha+n/p$ then $W\dot{K}_{p,\infty}^{\alpha}\subset\mathcal{M}%
_{q}^{r}$ (and the reverse) never could hold. This follows from an easy
scaling analysis of the space norms; in fact, the scaling of $\mathcal{M}%
_{q}^{r}$ is $-n/r$ and the one of $W\dot{K}_{p,\infty}^{0}$ is $-\alpha-n/p$.
\fin

In the next remark, we recall some inclusion and non-inclusion relations
involving Herz, weak-Herz, Besov and $bmo^{-1}$ spaces that can be found in
\cite{Tsutsui-1}.

\begin{remark}
\label{Rem Inclusoes tsutsui}

\begin{itemize}
\item[i)] For $1<p<\sigma<\infty$ and $0<\alpha<n\left(  1-1/p\right)  ,$ we
have
\[
W\dot{K}_{p,\infty}^{\alpha}\hookrightarrow\dot{B}_{\sigma,\infty}%
^{-(\alpha+n\left(  1/p-1/\sigma\right)  )}\text{, }\dot{K}_{p,\infty}%
^{\alpha}\hookrightarrow\dot{B}_{p,\infty}^{-\alpha}\text{ and }W\dot
{K}_{p,\sigma}^{0}\hookrightarrow\dot{B}_{\sigma,\infty}^{-n\left(
1/p-1/\sigma\right)  }.\text{ }%
\]

\item[ii)] For $1<p<\infty$ and $0\leq\alpha<n\left(  1-1/p\right)  ,$ we have
$W\dot{K}_{p,\infty}^{\alpha}\hookrightarrow\dot{B}_{\infty,\infty}%
^{-(\alpha+n/p)}.$

\item[iii)] For $0\leq\alpha<n,$ we have $\dot{K}_{\infty,\infty}^{\alpha
}\hookrightarrow\dot{B}_{\infty,\infty}^{-\alpha}.$

\item[iv)] For $1<p\leq\infty$ and $0\leq\alpha\leq n(1-1/p)$, we have
$W\dot{K}_{p,1}^{\alpha}\hookrightarrow\dot{B}_{\infty,\infty}^{-(\alpha
+n/p)}.$

\item[v)] We have $L^{1}=\dot{K}_{1,1}^{0}\hookrightarrow\dot{B}%
_{\infty,\infty}^{-n}.$ For $n<p\leq\infty$ and $0\leq\alpha<1-n/p$, the
inclusion $W\dot{K}_{p,\infty}^{\alpha}\hookrightarrow bmo^{-1}$ holds.

\item[vi)] For $1<p<\sigma<\infty$ and $-n\left(  1/p-1/\sigma\right)
<\alpha\leq0,$ $W\dot{K}_{p,\infty}^{\alpha}\hookrightarrow\dot{B}%
_{\sigma,\infty}^{-(\alpha+n\left(  1/p-1/\sigma\right)  )}$ does not hold.

\item[vii)] For $1<p<\infty$ and $-n/p<\alpha<0,$ $W\dot{K}_{p,\infty}%
^{\alpha}\hookrightarrow\dot{B}_{\infty,\infty}^{-(\alpha+n/p)}$ does not hold.
\end{itemize}
\end{remark}

\begin{remark}
\label{Remark. Relaciones entre BWK e Besov} Using the interpolation
properties of homogeneous Besov spaces and homogeneous Besov-weak-Herz spaces
(see Lemma \ref{Besov-weak-Herz-Interp} below) and item (ii) of Remark
\ref{Rem Inclusoes tsutsui}, for $1<p<\infty$ and $0\leq\alpha<n\left(
1-1/p\right)  $ we can obtain%
\begin{equation}
\dot{B}W\dot{K}_{p,\infty,r}^{\alpha,s}\hookrightarrow\dot{B}_{\infty
,r}^{s-(\alpha+n/p)}. \label{eq:Inclusao en Besov infty r s-(alpha+n/p)}%
\end{equation}
In particular, $\dot{B}W\dot{K}_{p,\infty,\infty}^{\alpha,\alpha+\frac{n}%
{p}-1}\hookrightarrow\dot{B}_{\infty,\infty}^{-1}$ and
\begin{equation}
\dot{B}W\dot{K}_{p,\infty,1}^{0,\frac{n}{p}}\hookrightarrow\dot{B}_{\infty
,1}^{0}\hookrightarrow L^{\infty}. \label{eq:imersao em Linfty}%
\end{equation}
Moreover, from Remark \ref{Rem Inclusoes tsutsui} (vi) and Lemma
\ref{Lem: Inclusiones Besov-herz con Herz} below, it follows that the
inclusion%
\[
\dot{B}W\dot{K}_{p,\infty,\infty}^{0,s}\hookrightarrow\dot{B}_{\sigma,\infty
}^{s-n\left(  1/p-1/\sigma\right)  }%
\]
does not hold for any $s\in\mathbb{R}$, $1<p<\infty$ and $1\leq\sigma<\infty$.
\end{remark}

\begin{remark}
Note that for $s-(\alpha+n/p)<0$ and $r>1$, or $s-(\alpha+n/p)\leq0$ and
$r=1$, the inclusion \eqref{eq:Inclusao en Besov infty r s-(alpha+n/p)}
implies that for $f\in\dot{B}W\dot{K}_{p,\infty,r}^{\alpha,s}$ the series
$\sum\limits_{j=-\infty}^{\infty}\Delta_{j}f$ converges in $\mathcal{S}%
^{\prime}$ to a representative of $f$ in $\mathcal{S}^{\prime}/\mathcal{P}$
(see e.g \cite{Lemarie}). So, in these cases the space $\dot{B}W\dot
{K}_{p,\infty,r}^{\alpha,s}$ can be regarded as a subspace of $\mathcal{S}%
^{\prime}$. Hereafter, we say that $f\in\mathcal{S}^{\prime}$ belongs to
$\dot{B}W\dot{K}_{p,\infty,r}^{\alpha,s}$ with $s-(\alpha+n/p)<0$ and $r>1$,
or $s-(\alpha+n/p)\leq0$ and $r=1$, if $f$ is the canonical representative of
the class in $\mathcal{S}^{\prime}/\mathcal{P}$, namely $f=\sum
\limits_{j=-\infty}^{\infty}\Delta_{j}f$ in $\mathcal{S}^{\prime}.$
\end{remark}

A multiplier theorem of H\"{o}rmander-Mihlin type will be needed in our
setting. This is the subject of the next lemma. In fact, the main part of the
proof has already been done in Lemma \ref{Teo:Operador-pseudo-Weak-Herz}.

\begin{lemma}
\label{Teo:Operador-pseudo-Besov-herz} Let $1<p<\infty,$ $1\leq q,r\leq
\infty,$ $-\frac{n}{p}<\alpha<n\left(  1-\frac{1}{p}\right)  $ and
$m,s\in\mathbb{R}$. Let $P\in C^{n}\left(  \mathbb{R}^{n}\backslash
\{0\}\right)  $ be a function such that $\left\vert \partial_{\xi}^{\beta
}P\left(  \xi\right)  \right\vert \leq C\left\vert \xi\right\vert ^{\left(
m-\left\vert \beta\right\vert \right)  }$ for all multi-index $\beta$
satisfying $\left\vert \beta\right\vert \leq n$. Then%
\[
\left\Vert P\left(  D\right)  f\right\Vert _{\dot{B}W\dot{K}_{p,q,r}%
^{\alpha,s-m}}\leq C\left\Vert f\right\Vert _{\dot{B}W\dot{K}_{p,q,r}%
^{\alpha,s}}.
\]

\end{lemma}

\textbf{Proof.} Note that for each $j\in\mathbb{Z}$ we have that $\left\vert
\xi\right\vert ^{m-\left\vert \beta\right\vert }\leq C2^{j\left(  m-\left\vert
\beta\right\vert \right)  }$ for all $\xi\in\tilde{D}_{j}$, and therefore
$\left\vert \partial_{\xi}^{\beta}P\left(  \xi\right)  \right\vert \leq
C2^{j\left(  m-\left\vert \beta\right\vert \right)  }.$ On the other hand,
since $\mbox{supp}\widehat{\Delta_{j}f}\subset D_{j}$ we can use Lemma
\ref{Teo:Operador-pseudo-Weak-Herz} in order to get%

\begin{equation}
\left\Vert \Delta_{j}\left(  P\left(  D\right)  f\right)  \right\Vert
_{W\dot{K}_{p,q}^{\alpha}}=\left\Vert P\left(  D\right)  \left(  \Delta
_{j}f\right)  \right\Vert _{W\dot{K}_{p,q}^{\alpha}}\leq C2^{jm}\left\Vert
\Delta_{j}f\right\Vert _{W\dot{K}_{p,q}^{\alpha}}. \label{aux-pseudo-1}%
\end{equation}
Now the result follows by multiplying (\ref{aux-pseudo-1}) by $2^{j\left(
s-m\right)  }$ and after taking the $l^{r}$-norm.\fin\bigskip

In what follows we present some inclusions involving Sobolev-weak-Herz and
Besov-weak-Herz spaces.

\begin{lemma}
\label{Lem: Inclusiones Besov-herz con Herz} Let $s\in\mathbb{R},$
$1<p<\infty,$ $1\leq q\leq\infty$ and $-\frac{n}{p}<\alpha<n\left(  1-\frac
{1}{p}\right)  $. We have the following continuous inclusions%

\begin{align}
\dot{B}W\dot{K}_{p,q,1}^{\alpha,0}  &  \subset W\dot{K}_{p,q}^{\alpha}%
\subset\dot{B}W\dot{K}_{p,q,\infty}^{\alpha,0}%
\label{eq:Inclusion De BK com s=00003D00003D00003D00003D0}\\
\dot{B}W\dot{K}_{p,q,1}^{\alpha,s}  &  \subset W\dot{K}_{p,q}^{\alpha
,s}\subset\dot{B}W\dot{K}_{p,q,\infty}^{\alpha,s}.
\label{eq:Inclusion de BK con s}%
\end{align}

\end{lemma}

\textbf{Proof. }For $f\in\dot{B}W\dot{K}_{p,q,1}^{\alpha,0}$, we can employ
the decomposition $f=\sum\limits_{j\in\mathbb{Z}}\Delta_{j}f$ in order to estimate%

\[
\left\Vert f\right\Vert _{L^{p,\infty}\left(  A_{k}\right)  }\leq
\sum\limits_{j\in\mathbb{Z}}\left\Vert \Delta_{j}f\right\Vert _{L^{p,\infty
}\left(  A_{k}\right)  }.
\]
Thus, using Minkowski inequality, we arrive at (with usual modification in the
case $q=\infty$)%

\begin{align*}
\left\Vert f\right\Vert _{W\dot{K}_{p,q}^{\alpha}}  &  \leq\left[
\sum\limits_{k\in\mathbb{Z}}2^{k\alpha q}\left\Vert f\right\Vert
_{L^{p,\infty}\left(  A_{k}\right)  }^{q}\right]  ^{\frac{1}{q}}\leq\left[
\sum\limits_{k\in\mathbb{Z}}\left(  \sum\limits_{j\in\mathbb{Z}}2^{k\alpha
}\left\Vert \Delta_{j}f\right\Vert _{L^{p,\infty}\left(  A_{k}\right)
}\right)  ^{q}\right]  ^{\frac{1}{q}}\\
&  \leq\sum\limits_{j\in\mathbb{Z}}\left(  \sum\limits_{k\in\mathbb{Z}%
}2^{k\alpha q}\left\Vert \Delta_{j}f\right\Vert _{L^{p,\infty}\left(
A_{k}\right)  }^{q}\right)  ^{\frac{1}{q}}=\sum\limits_{j\in\mathbb{Z}%
}\left\Vert \Delta_{j}f\right\Vert _{W\dot{K}_{p,q}^{\alpha}}\\
&  =\left\Vert f\right\Vert _{\dot{B}W\dot{K}_{p,q,1}^{\alpha,0}},
\end{align*}
which implies the first inclusion in
(\ref{eq:Inclusion De BK com s=00003D00003D00003D00003D0}). Now, let $f\in
W\dot{K}_{p,q}^{\alpha}$ and note that in fact we have that $f\in
\mathcal{\mathcal{S}}^{\prime}/\mathcal{P}$. Moreover, using Lemma
\ref{lem:[Comvolucao]} we get%

\[
\left\Vert f\right\Vert _{\dot{B}W\dot{K}_{p,q,\infty}^{\alpha,0}%
}=\mathop{\sup}\limits_{j\in\mathbb{Z}}\left\Vert \Delta_{j}f\right\Vert
_{W\dot{K}_{p,q}^{\alpha}}\leq C\mathop{\sup}\limits_{j\in\mathbb{Z}%
}\left\Vert f\right\Vert _{W\dot{K}_{p,q}^{\alpha}}=C\left\Vert f\right\Vert
_{W\dot{K}_{p,q}^{\alpha}},
\]
and then the second inclusion in
(\ref{eq:Inclusion De BK com s=00003D00003D00003D00003D0}) holds.

For (\ref{eq:Inclusion de BK con s}), we can use Lemma
\ref{Teo:Operador-pseudo-Weak-Herz} in order to estimate%

\begin{align*}
\left\Vert f\right\Vert _{W\dot{K}_{p,q}^{\alpha,s}}  &  =\left\Vert
I^{s}f\right\Vert _{W\dot{K}_{p,q}^{\alpha}}\leq\left\Vert I^{s}f\right\Vert
_{\dot{B}W\dot{K}_{p,q,1}^{\alpha,0}}=\sum\limits_{j\in\mathbb{Z}}\left\Vert
\Delta_{j}I^{s}f\right\Vert _{W\dot{K}_{p,q}^{\alpha}}\\
&  \leq C\sum\limits_{j\in\mathbb{Z}}2^{js}\left\Vert \Delta_{j}f\right\Vert
_{W\dot{K}_{p,q}^{\alpha}}=C\left\Vert f\right\Vert _{\dot{B}W\dot{K}%
_{p,q,1}^{\alpha,s}}\text{.}%
\end{align*}
Moreover, Lemma \ref{Teo:Operador-pseudo-Weak-Herz} also can be used to obtain%

\begin{align*}
\left\Vert f\right\Vert _{\dot{B}W\dot{K}_{p,q,\infty}^{\alpha,s}}  &
=\mathop{\sup}\limits_{j\in\mathbb{Z}}2^{js}\left\Vert \Delta_{j}f\right\Vert
_{W\dot{K}_{p,q}^{\alpha}}=\mathop{\sup}\limits_{j\in\mathbb{Z}}%
2^{js}\left\Vert I^{-s}\Delta_{j}I^{s}f\right\Vert _{W\dot{K}_{p,q}^{\alpha}%
}\\
&  \leq C\mathop{\sup}\limits_{j\in\mathbb{Z}}\left\Vert \Delta_{j}%
I^{s}f\right\Vert _{W\dot{K}_{p,q}^{\alpha}}\leq C\mathop{\sup}\limits_{j\in
\mathbb{Z}}\left\Vert I^{s}f\right\Vert _{W\dot{K}_{p,q}^{\alpha}}\\
&  =C\left\Vert I^{s}f\right\Vert _{W\dot{K}_{p,q}^{\alpha}}=C\left\Vert
f\right\Vert _{W\dot{K}_{p,q}^{\alpha,s}},
\end{align*}
for all $f\in W\dot{K}_{p,q}^{\alpha,s}$, as required.

\fin

Now we present an embedding theorem of Sobolev type.

\begin{lemma}
\label{Sobolev-Lemma-1}Let $s\in\mathbb{R}$, $1<p<\infty$, $1\leq
q,r\leq\infty,$ $p\leq p_{1}<\infty$, $1<p_{2}\leq p_{1}$ and $-\frac{n}%
{p}<\alpha<n\left(  1+\frac{1}{p_{1}}-\frac{1}{p_{2}}-\frac{1}{p}\right)  $. Then%

\begin{equation}
\left\Vert f\right\Vert _{\dot{B}W\dot{K}_{p,q,r}^{\alpha,s}}\leq C\left\Vert
f\right\Vert _{\dot{B}W\dot{K}_{p_{2},q,r}^{\alpha+n\left(  \frac{1}{p}%
-\frac{1}{p_{1}}\right)  ,s+n\left(  \frac{1}{p_{2}}-\frac{1}{p_{1}}\right)
}}. \label{eq:Inclucao geral besov-WH}%
\end{equation}

In particular, for $\frac{n}{2}<p<\infty$ and $0\le\alpha<\min\left\{
1-\frac{n}{2p},\frac{n}{2p}\right\}  $, it follows that%

\begin{equation}
\left\Vert f\right\Vert _{\dot{B}W\dot{K}_{2p,q,r}^{\alpha,s}}\leq C\left\Vert
f\right\Vert _{\dot{B}W\dot{K}_{p,q,r}^{2\alpha,\alpha+s+\frac{n}{2p}}}.
\label{eq:Inclucao dobrando alpha}%
\end{equation}

\end{lemma}

\textbf{Proof.} Using H\"{o}lder inequality, it follows that%

\[
\left\Vert \Delta_{j}f\right\Vert _{W\dot{K}_{p,q,r}^{\alpha}}\le C\left\Vert
\Delta_{j}f\right\Vert _{W\dot{K}_{p_{1},q,r}^{\alpha+n\left(  \frac{1}%
{p}-\frac{1}{p_{1}}\right)  }}.
\]

Also, we have that $\varphi_{j}\hat{f}=\tilde{\varphi}_{j}\varphi_{j}\hat{f}$,
that is, $\Delta_{j}f=\left(  \tilde{\varphi}_{j}\right)  \check{}\ast
\Delta_{j}f.$ So, using Lemma \ref{lem:[Comvolucao]} we get%

\begin{align*}
\left\Vert \Delta_{j}f\right\Vert _{W\dot{K}_{p_{1},q,r}^{\alpha+n\left(
\frac{1}{p}-\frac{1}{p_{1}}\right)  }}  &  =\left\Vert \left(  \tilde{\varphi
}_{j}\right)  \check{}\ast\Delta_{j}f\right\Vert _{W\dot{K}_{p_{1}%
,q,r}^{\alpha+n\left(  \frac{1}{p}-\frac{1}{p_{1}}\right)  }}\\
&  \leq C\max\left\{  \left\Vert \left(  \tilde{\varphi}_{j}\right)  \check
{}\right\Vert _{L^{p\ast}},\left\Vert \left\vert \cdot\right\vert ^{\frac
{n}{p\ast}}\left(  \tilde{\varphi}_{j}\right)  \check{}\right\Vert
_{L^{\infty}}\right\}  \left\Vert \Delta_{j}f\right\Vert _{W\dot{K}%
_{p_{2},q,r}^{\alpha+n\left(  \frac{1}{p}-\frac{1}{p_{1}}\right)  }},
\end{align*}
where $1+\frac{1}{p_{1}}=\frac{1}{p\ast}+\frac{1}{p_{2}}$. It is easy to check
that $\max\left\{  \left\Vert \left(  \tilde{\varphi}_{j}\right)  \check
{}\right\Vert _{L^{p\ast}},\left\Vert \left\vert \cdot\right\vert ^{n/p\ast
}\left(  \tilde{\varphi}_{j}\right)  \check{}\right\Vert _{L^{\infty}%
}\right\}  \leq C2^{jn\left(  \frac{1}{p_{2}}-\frac{1}{p_{1}}\right)  }$, and then%

\[
\left\Vert \Delta_{j}f\right\Vert _{W\dot{K}_{p_{1},q,r}^{\alpha+n\left(
\frac{1}{p}-\frac{1}{p_{1}}\right)  }}\leq C2^{jn\left(  \frac{1}{p_{2}}%
-\frac{1}{p_{1}}\right)  }\left\Vert \Delta_{j}f\right\Vert _{W\dot{K}%
_{p_{2},q,r}^{\alpha+n\left(  \frac{1}{p}-\frac{1}{p_{1}}\right)  }},
\]
which gives \eqref{eq:Inclucao geral besov-WH}. We conclude the proof by
noting that for $0\leq\alpha<n/2p$ there exists $p_{1}$ such that $p_{1}%
\geq2p$ and $\alpha=n\left(  \frac{1}{p}-\frac{1}{p_{1}}\right)  $. Moreover,
$\alpha<n+\frac{n}{p_{1}}-\frac{1}{p}-\frac{n}{2p}$ because $\alpha<1-\frac
{n}{2p}\leq\frac{n}{2}-\frac{n}{2p}$. So, \eqref{eq:Inclucao dobrando alpha}
follows from \eqref{eq:Inclucao geral besov-WH} by choosing this value of
$p_{1}$.

\fin

We finish this section with a result that provides a characterization of
homogeneous Besov-weak-Herz spaces as interpolation of two homogeneous
Sobolev-weak-Herz ones.

\begin{lemma}
\label{Besov-weak-Herz-Interp} Let $s_{0},s_{1},s\in\mathbb{R}$, $1<p<\infty$,
$1\leq q,r\leq\infty$ and $-\frac{n}{p}<\alpha<n\left(  1-\frac{1}{p}\right)
.$ If $s_{0}\neq s_{1}$ and $s=(1-\theta)s_{0}+\theta s_{1}$ with $\theta
\in\left(  0,1\right)  $, then%

\[
\left(  W\dot{K}_{p,q}^{\alpha,s_{0}},W\dot{K}_{p,q}^{\alpha,s_{1}}\right)
_{\theta,r}=\dot{B}W\dot{K}_{p,q,r}^{\alpha,s}.
\]

\end{lemma}

\textbf{Proof.} Let $f=f_{0}+f_{1}$ with $f_{i}\in W\dot{K}_{p,q}%
^{\alpha,s_{i}}$ $i=0,1$. By using the Lemma
\ref{Teo:Operador-pseudo-Weak-Herz} we get%

\begin{align}
\left\Vert \Delta_{j}f\right\Vert _{W\dot{K}_{p,q}^{\alpha}}  &
\leq\left\Vert \Delta_{j}f_{0}\right\Vert _{W\dot{K}_{p,q}^{\alpha}%
}+\left\Vert \Delta_{j}f_{1}\right\Vert _{W\dot{K}_{p,q}^{\alpha}}\nonumber\\
&  \leq C\left(  2^{-s_{0}j}\left\Vert I^{s_{0}}\Delta_{j}f_{0}\right\Vert
_{W\dot{K}_{p,q}^{\alpha}}+2^{-s_{1}j}\left\Vert I^{s_{1}}\Delta_{j}%
f_{1}\right\Vert _{W\dot{K}_{p,q}^{\alpha}}\right) \nonumber\\
&  \leq C\left(  2^{-s_{0}j}\left\Vert I^{s_{0}}f_{0}\right\Vert _{W\dot
{K}_{p,q}^{\alpha}}+2^{-s_{1}j}\left\Vert I^{s_{1}}f_{1}\right\Vert _{W\dot
{K}_{p,q}^{\alpha}}\right) \nonumber\\
&  \leq C2^{-s_{0}j}\left(  \left\Vert f_{0}\right\Vert _{W\dot{K}%
_{p,q}^{\alpha,s_{0}}}+2^{\left(  s_{0}-s_{1}\right)  j}\left\Vert
f_{1}\right\Vert _{W\dot{K}_{p,q}^{\alpha,s_{1}}}\right)  .
\label{aux-interp-1000}%
\end{align}
It follows from (\ref{aux-interp-1000}) that%

\[
\left\Vert \Delta_{j}f\right\Vert _{W\dot{K}_{p,q}^{\alpha}}\leq C2^{-s_{0}%
j}K\left(  2^{\left(  s_{0}-s_{1}\right)  j},f,W\dot{K}_{p,q}^{\alpha,s_{0}%
},W\dot{K}_{p,q}^{\alpha,s_{1}}\right)  .
\]
Noting that $s-s_{0}=-\theta\left(  s_{0}-s_{1}\right)  $ and multiplying by
$2^{js}$ the previous inequality, we arrive at%

\[
2^{sj}\left\Vert \Delta_{j}f\right\Vert _{W\dot{K}_{p,q}^{\alpha}}\leq
C\left(  2^{\left(  s_{0}-s_{1}\right)  j}\right)  ^{-\theta}K\left(
2^{\left(  s_{0}-s_{1}\right)  j},f,W\dot{K}_{p,q}^{\alpha,s_{0}},W\dot
{K}_{p,q}^{\alpha,s_{1}}\right)  ,
\]
and then (see \cite[Lemma 3.1.3]{BL}) we can conclude that%

\[
\left\Vert f\right\Vert _{\dot{B}W\dot{K}_{p,q,r}^{\alpha,s}}\le C\left\Vert
f\right\Vert _{\left(  W\dot{K}_{p,q}^{\alpha,s_{0}},W\dot{K}_{p,q}%
^{\alpha,s_{1}}\right)  _{\theta,r}}.
\]

In order to prove the reverse inequality, note that by using again Lemma
\ref{Teo:Operador-pseudo-Weak-Herz} we have%

\begin{align*}
2^{\left(  s-s_{0}\right)  j}J\left(  2^{\left(  s_{0}-s_{1}\right)  j}%
,\Delta_{j}f,W\dot{K}_{p,q}^{\alpha,s_{0}},W\dot{K}_{p,q}^{\alpha,s_{1}%
}\right)   &  =2^{\left(  s-s_{0}\right)  j}\mbox{max}\left\{  \left\Vert
\Delta_{j}f\right\Vert _{W\dot{K}_{p,q}^{\alpha,s_{0}}},2^{\left(  s_{0}%
-s_{1}\right)  j}\left\Vert \Delta_{j}f\right\Vert _{W\dot{K}_{p,q}%
^{\alpha,s_{1}}}\right\} \\
&  \leq2^{\left(  s-s_{0}\right)  j}\mbox{max}\left\{  2^{s_{0}j}\left\Vert
\Delta_{j}f\right\Vert _{W\dot{K}_{p,q}^{\alpha}},2^{s_{0}j}\left\Vert
\Delta_{j}f\right\Vert _{W\dot{K}_{p,q}^{\alpha}}\right\} \\
&  \leq2^{sj}\mbox{max}\left\{  \left\Vert \Delta_{j}f\right\Vert _{W\dot
{K}_{p,q}^{\alpha}},\left\Vert \Delta_{j}f\right\Vert _{W\dot{K}_{p,q}%
^{\alpha}}\right\} \\
&  =2^{sj}\left\Vert \Delta_{j}f\right\Vert _{W\dot{K}_{p,q}^{\alpha}}.
\end{align*}
Now the Equivalence Theorem (see \cite[Lemma 3.2.3]{BL}) leads us to%

\[
\left\Vert f\right\Vert _{\left(  W\dot{K}_{p,q}^{\alpha,s_{0}},W\dot{K}%
_{p,q}^{\alpha,s_{1}}\right)  _{\theta,r}}\le C\left\Vert f\right\Vert
_{\dot{B}W\dot{K}_{p,q,r}^{\alpha,s}}.
\]

The remainder of the proof is to show that in fact $f\in\dot{B}W\dot
{K}_{p,q,r}^{\alpha,s}$ implies that $f\in W\dot{K}_{p,q}^{\alpha,s_{0}}%
+W\dot{K}_{p,q}^{\alpha,s_{1}}$. Suppose that $s_{0}>s_{1}$ (without loss of
generality). Using the decomposition $f=\sum\limits_{j<0}\Delta_{j}%
f+\sum\limits_{j\geq0}\Delta_{j}f=f_{0}+f_{1}$ and Lemma
\ref{Teo:Operador-pseudo-Weak-Herz}, we obtain%

\begin{align*}
\left\Vert f_{0}\right\Vert _{W\dot{K}_{p,q}^{\alpha,s_{0}}}  &  \leq
\sum\limits_{j<0}\left\Vert \Delta_{j}f\right\Vert _{W\dot{K}_{p,q}%
^{\alpha,s_{0}}}\leq\sum\limits_{j<0}2^{j(s_{0}-s)}2^{js}\left\Vert \Delta
_{j}f\right\Vert _{W\dot{K}_{p,q}^{\alpha}}\\
&  \leq C\left(  \sum\limits_{j<0}2^{j(s_{0}-s)r^{\prime}}\right)  ^{\frac
{1}{r^{\prime}}}\left(  \sum\limits_{j<0}2^{jsr}\left\Vert \Delta
_{j}f\right\Vert _{W\dot{K}_{p,q}^{\alpha}}^{r}\right)  ^{\frac{1}{r}}\\
&  \leq C\left\Vert f\right\Vert _{\dot{B}W\dot{K}_{p,q,r}^{\alpha,s}}.
\end{align*}
Similarly, one has%

\begin{align*}
\left\Vert f_{1}\right\Vert _{W\dot{K}_{p,q}^{\alpha,s_{1}}}  &  \leq
\sum\limits_{j\geq0}\left\Vert \Delta_{j}f\right\Vert _{W\dot{K}_{p,q}%
^{\alpha,s_{1}}}\leq\sum\limits_{j\geq0}2^{j(s_{1}-s)}2^{js}\left\Vert
\Delta_{j}f\right\Vert _{W\dot{K}_{p,q}^{\alpha}}\\
&  \leq C\left(  \sum\limits_{j\geq0}2^{j(s_{1}-s)r^{\prime}}\right)
^{1/r^{\prime}}\left(  \sum\limits_{j\geq0}2^{jsr}\left\Vert \Delta
_{j}f\right\Vert _{W\dot{K}_{p,q}^{\alpha}}^{r}\right)  ^{1/r}\\
&  \leq C\left\Vert f\right\Vert _{\dot{B}W\dot{K}_{p,q,r}^{\alpha,s}}%
\end{align*}
and then we are done.\fin

\bigskip

\section{Proof of Theorem \ref{Teo-main}}

In the previous sections, we have derived key properties about homogeneous
Besov-weak-Herz spaces. With these results in hands, we prove Theorem
\ref{Teo-main} in the present section.

\subsection{Heat kernel estimates}

We start by providing estimates for the heat semigroup $\{G(t)\}_{t\geq0}$ in
Besov-weak-Herz spaces. Recall that in the whole space $\mathbb{R}^{n}$ this
semigroup can be defined as $G(t)f=\left(  \exp\left(  -t\left\vert
\xi\right\vert ^{2}\right)  \hat{f}\right)  \check{}$ $\,$ for all
$f\in\mathcal{S}^{\prime}$ and $t\geq0.$

\begin{lemma}
\label{Heat-Kernel-estimate} Let $s,\sigma\in\mathbb{R}$, $s\leq\sigma$,
$1<p<\infty,$ $1\leq q,r\leq\infty$ and $-\frac{n}{p}<\alpha<n\left(
1-\frac{1}{p}\right)  $. Then, there is $C>0$ (independent of $f$) such that
\begin{equation}
\left\Vert G\left(  t\right)  f\right\Vert _{\dot{B}W\dot{K}_{p,q,r}%
^{\alpha,\sigma}}\leq Ct^{\left(  s-\sigma\right)  /2}\left\Vert f\right\Vert
_{\dot{B}W\dot{K}_{p,q,r}^{\alpha,s}},
\label{eq:Heat-Kernel estimate  s leq sigma}%
\end{equation}
for all $t>0.$ Moreover, if $s<\sigma$, then we have the estimate%

\begin{equation}
\left\Vert G\left(  t\right)  f\right\Vert _{\dot{B}W\dot{K}_{p,q,1}%
^{\alpha,\sigma}}\leq Ct^{\left(  s-\sigma\right)  /2}\left\Vert f\right\Vert
_{\dot{B}W\dot{K}_{p,q,\infty}^{\alpha,s}},
\label{eq:Heat-Kernel estimate  s menor sigma}%
\end{equation}
for all $t>0.$
\end{lemma}

\textbf{Proof.} Firstly, observe that for each multi-index $\beta$ there is a
polynomial $p_{\beta}(\cdot)$ of degree $\left\vert \beta\right\vert $ such that%

\[
\partial_{\xi}^{\beta}\left(  \exp(-t\left\vert \xi\right\vert ^{2})\right)
=t^{\left\vert \beta\right\vert /2}p_{\beta}(\sqrt{t}\xi)\exp\left(
-t\left\vert \xi\right\vert ^{2}\right)  .
\]
Therefore, for some $C>0$ it follows that%

\[
\left\vert \partial_{\xi}^{\beta}\left(  \exp(-t\left\vert \xi\right\vert
^{2})\right)  \right\vert \leq Ct^{-m/2}\left\vert \xi\right\vert
^{-m-\left\vert \beta\right\vert }.
\]
By employing Lemma \ref{Teo:Operador-pseudo-Besov-herz}, we obtain%

\[
\left\Vert G\left(  t\right)  f\right\Vert _{\dot{B}W\dot{K}_{p,q,r}%
^{\alpha,s-m}}\leq Ct^{-m/2}\left\Vert f\right\Vert _{\dot{B}W\dot{K}%
_{p,q,r}^{\alpha,s}}.
\]
Taking now $m=s-\sigma$ we arrive at the inequality \eqref{eq:Heat-Kernel estimate  s leq sigma}.

Next we turn to \eqref{eq:Heat-Kernel estimate  s menor sigma} and let
$s<\sigma.$ From \eqref{eq:Heat-Kernel estimate  s leq sigma} with $r=\infty$
we get%

\[
\left\Vert G\left(  t\right)  f\right\Vert _{\dot{B}W\dot{K}_{p,q,\infty
}^{\alpha,2\sigma-s}}\leq Ct^{s-\sigma}\left\Vert f\right\Vert _{\dot{B}%
W\dot{K}_{p,q,\infty}^{\alpha,s}}%
\]
and%

\[
\left\Vert G\left(  t\right)  f\right\Vert _{\dot{B}W\dot{K}_{p,q,\infty
}^{\alpha,s}}\leq C\left\Vert f\right\Vert _{\dot{B}W\dot{K}_{p,q,\infty
}^{\alpha,s}}.
\]
By using Lemma \ref{Besov-weak-Herz-Interp} and the Reiteration Theorem (see
\cite[Theorem 3.5.3 and its remark]{BL}) we conclude that%

\[
G(t):\,\dot{B}W\dot{K}_{p,q,\infty}^{\alpha,s}\rightarrow\left(  \dot{B}%
W\dot{K}_{p,q,\infty}^{\alpha,2\sigma-s},\dot{B}W\dot{K}_{p,q,\infty}%
^{\alpha,s}\right)  _{\frac{1}{2},1}=\dot{B}W\dot{K}_{p,q,1}^{\alpha,\sigma},
\]
with $\left\Vert G\left(  t\right)  \right\Vert _{\dot{B}W\dot{K}_{p,q,\infty
}^{\alpha,s}\rightarrow\dot{B}W\dot{K}_{p,q,1}^{\alpha,\sigma}}\leq
Ct^{\left(  s-\sigma\right)  /2}$, which gives \eqref{eq:Heat-Kernel estimate  s menor sigma}.

\fin

\subsection{Bilinear estimate}

\bigskip Let us define the space $X$ as%

\[
X=\left\{  u:\,(0,\infty)\rightarrow\dot{B}W\dot{K}_{p,q,\infty}
^{\alpha,\alpha+\frac{n}{p}-1}\cap W\dot{K}_{2p,2q}^{\alpha}%
\,\mbox{with \ensuremath{\nabla\cdot u=0}}\,\mbox{such that}\,\left\Vert
u\right\Vert _{X}<\infty\right\}  ,
\]
where%

\begin{equation}
\left\Vert u\right\Vert _{X}:=\left\Vert u\right\Vert _{L^{\infty}
((0,\infty);\dot{B}W\dot{K}_{p,q,\infty} ^{\alpha,\alpha+\frac{n}{p}-1}%
)}+\mathop{\sup}\limits_{t>0}t^{\frac{1}{2}-\left(  \frac{\alpha}{2}+\frac
{n}{4p}\right)  }\left\Vert u\right\Vert _{W\dot{K}_{2p,2q}^{\alpha}}.
\label{aux-norm-1}%
\end{equation}

We are going to prove the bilinear estimate%

\begin{equation}
\left\Vert B(u,v)\right\Vert _{X}\leq K\left\Vert u\right\Vert _{X}\left\Vert
v\right\Vert _{X}. \label{eq:Estimativa termo naolinear em X}%
\end{equation}

We start by estimating the second part of the norm (\ref{aux-norm-1}). For
that matter, we use
\eqref{eq:Inclusion De BK com s=00003D00003D00003D00003D0},
\eqref{eq:Inclucao dobrando alpha},
\eqref{eq:Heat-Kernel estimate  s menor sigma} and Lemma
\ref{Teo:Operador-pseudo-Besov-herz} in order to get%

\begin{align*}
\left\Vert B(u,v)(t)\right\Vert _{W\dot{K}_{2p,2q}^{\alpha}}  &
\leq\left\Vert B(u,v)(t)\right\Vert _{\dot{B}W\dot{K}_{2p,2q,1}^{\alpha,0}%
}\leq\left\Vert B(u,v)(t)\right\Vert _{\dot{B}W\dot{K}_{p,2q,1}^{2\alpha
,\alpha+\frac{n}{2p}}}\\
&  \leq C\int_{0}^{t}\left\Vert G\left(  t-\tau\right)  \mathbb{P}%
\mbox{{\rm div}}\left(  u\otimes v\right)  \right\Vert _{\dot{B}W\dot
{K}_{p,2q,1}^{2\alpha,\alpha+\frac{n}{2p}}}d\tau\\
&  \leq C\int_{0}^{t}\left(  t-\tau\right)  ^{-\frac{1}{2}-\left(
\frac{\alpha}{2}+\frac{n}{4p}\right)  }\left\Vert \mathbb{P}%
\mbox{{\rm div}}\left(  u\otimes v\right)  \right\Vert _{\dot{B}W\dot
{K}_{p,2q,\infty}^{2\alpha,-1}}d\tau\\
&  \leq C\int_{0}^{t}\left(  t-\tau\right)  ^{-\frac{1}{2}-\left(
\frac{\alpha}{2}+\frac{n}{4p}\right)  }\left\Vert u\otimes v\right\Vert
_{\dot{B}W\dot{K}_{p,2q,\infty}^{2\alpha,0}}d\tau\\
&  \leq C\int_{0}^{t}\left(  t-\tau\right)  ^{-\frac{1}{2}-\left(
\frac{\alpha}{2}+\frac{n}{4p}\right)  }\left\Vert u\otimes v\right\Vert
_{W\dot{K}_{p,q}^{2\alpha}}d\tau\\
&  \leq C\int_{0}^{t}\left(  t-\tau\right)  ^{-\frac{1}{2}-\left(
\frac{\alpha}{2}+\frac{n}{4p}\right)  }\left\Vert u\right\Vert _{W\dot
{K}_{2p,2q}^{\alpha}}\left\Vert v\right\Vert _{W\dot{K}_{2p,2q}^{\alpha}}%
d\tau\\
&  \leq C\int_{0}^{t}\left(  t-\tau\right)  ^{-\frac{1}{2}-\left(
\frac{\alpha}{2}+\frac{n}{4p}\right)  }\tau^{-2\left(  \frac{1}{2}-\left(
\frac{\alpha}{2}+\frac{n}{4p}\right)  \right)  }d\tau\left\Vert u\right\Vert
_{X}\left\Vert v\right\Vert _{X}\\
&  \leq Ct^{^{-\frac{1}{2}+\left(  \frac{\alpha}{2}+\frac{n}{4p}\right)  }%
}\mathcal{B}\left(  \alpha+\frac{n}{2p},\frac{1}{2}-\left(  \frac{\alpha}%
{2}+\frac{n}{4p}\right)  \right)  \left\Vert u\right\Vert _{X}\left\Vert
v\right\Vert _{X},
\end{align*}
where $\mathcal{B}\left(  \cdot,\cdot\right)  $ denotes the beta function. The
previous estimate leads us to%

\begin{equation}
\mathop{\sup}\limits_{t>0}t^{\frac{1}{2}-\left(  \frac{\alpha}{2}+\frac{n}%
{4p}\right)  }\left\Vert B(u,v)(t)\right\Vert _{W\dot{K}_{2p,2q}^{\alpha}}\leq
C\left\Vert u\right\Vert _{X}\left\Vert v\right\Vert _{X}.
\label{eq:Estimativa termo naolinear-1}%
\end{equation}
Moreover, for the first part of the norm (\ref{aux-norm-1}), we have%

\begin{align*}
\left\Vert B(u,v)(t)\right\Vert _{\dot{B}W\dot{K}_{p,q,\infty}^{\alpha
,\alpha+\frac{n}{p}-1}}  &  \leq\int_{0}^{t}\left\Vert G\left(  t-\tau\right)
\mathbb{P}\mbox{{\rm div}}\left[  u\otimes v\right]  \right\Vert _{\dot
{B}W\dot{K}_{p,q,\infty}^{\alpha,\alpha+\frac{n}{p}-1}}d\tau\\
&  \leq C\int_{0}^{t}\left\Vert G\left(  t-\tau\right)  \mathbb{P}%
\mbox{{\rm div}}\left[  u\otimes v\right]  \right\Vert _{\dot{B}W\dot
{K}_{p,q,\infty}^{2\alpha,2\alpha+\frac{n}{p}-1}}d\tau\\
&  \leq C\int_{0}^{t}\left(  t-\tau\right)  ^{-\left(  \alpha+\frac{n}%
{2p}\right)  }\left\Vert \mathbb{P}\mbox{{\rm div}}\left[  u\otimes v\right]
\right\Vert _{\dot{B}W\dot{K}_{p,q,\infty}^{2\alpha,-1}}d\tau\\
&  \leq C\int_{0}^{t}\left(  t-\tau\right)  ^{-\left(  \alpha+\frac{n}%
{2p}\right)  }\left\Vert u\otimes v\right\Vert _{\dot{B}W\dot{K}_{p,q,\infty
}^{2\alpha,0}}d\tau\\
&  \leq C\int_{0}^{t}\left(  t-\tau\right)  ^{-\left(  \alpha+\frac{n}%
{2p}\right)  }\left\Vert u\otimes v\right\Vert _{W\dot{K}_{p,q}^{2\alpha}%
}d\tau\\
&  \leq C\left\Vert u\right\Vert _{X}\left\Vert v\right\Vert _{X}\int_{0}%
^{t}\left(  t-\tau\right)  ^{-\left(  \alpha+\frac{n}{2p}\right)  }%
\tau^{-2\left(  \frac{1}{2}-\left(  \frac{\alpha}{2}+\frac{n}{4p}\right)
\right)  }d\tau\\
&  \leq C\mathcal{B}\left(  \alpha+\frac{n}{2p},1-\left(  \alpha+\frac{n}%
{2p}\right)  \right)  \left\Vert u\right\Vert _{X}\left\Vert v\right\Vert
_{X}.
\end{align*}
In other words, we have obtained the estimate%

\begin{equation}
\left\Vert B(u,v)\right\Vert _{L^{\infty}(\left(  0,\infty\right)  ;\dot
{B}W\dot{K}_{p,q,\infty}^{\alpha,\alpha+\frac{n}{p}-1})}\leq C\left\Vert
u\right\Vert _{X}\left\Vert v\right\Vert _{X}.
\label{eq:Estimativa termo naolinear-2}%
\end{equation}
Finally, notice that the estimates \eqref{eq:Estimativa termo naolinear-1} and
\eqref{eq:Estimativa termo naolinear-2} together give
(\ref{eq:Estimativa termo naolinear em X}).

\subsection{Proof of Theorem \ref{Teo-main}}

\textbf{Existence and Uniqueness.} For $\epsilon>0$ (to be chosen later) let
$\bar{B}(0,\epsilon)$ denote the closed ball in $X$ and define the operator
$\Psi:\,\bar{B}(0,2\epsilon)\rightarrow\bar{B}(0,2\epsilon)$ as%

\[
\Psi(u)=G(t)u_{0}+B(u,u).
\]

First, note that by using
\eqref{eq:Inclusion De BK com s=00003D00003D00003D00003D0},
\eqref{eq:Heat-Kernel estimate  s menor sigma} , $\alpha+\frac{n}{2p}-1<0$ and
\eqref{eq:Inclucao geral besov-WH} it follows that%

\begin{align}
\mathop{\sup}\limits_{t>0}t^{\frac{1}{2}-\left(  \frac{\alpha}{2}+\frac{n}%
{4p}\right)  }\left\Vert G(t)u_{0}\right\Vert _{W\dot{K}_{2p,2q}^{\alpha}} &
\leq C\mathop{\sup}\limits_{t>0}t^{\frac{1}{2}-\left(  \frac{\alpha}{2}%
+\frac{n}{4p}\right)  }\left\Vert G(t)u_{0}\right\Vert _{\dot{B}W\dot
{K}_{2p,2q,1}^{\alpha,0}}\nonumber\\
&  \leq C\left\Vert u_{0}\right\Vert _{\dot{B}W\dot{K}_{2p,2q,\infty}%
^{\alpha,\alpha+\frac{n}{2p}-1}}\leq C\left\Vert u_{0}\right\Vert _{\dot
{B}W\dot{K}_{p,q,\infty}^{\alpha,\alpha+\frac{n}{p}-1}}.
\label{eq:Estimativa termo linear-1}%
\end{align}
Moreover, using \eqref{eq:Heat-Kernel estimate  s leq sigma} we obtain%

\[
\left\Vert G(t)u_{0}\right\Vert _{\dot{B}W\dot{K}_{p,q,\infty}^{\alpha
,\alpha+\frac{n}{p}-1}}\leq C\left\Vert u_{0}\right\Vert _{\dot{B}W\dot
{K}_{p,q,\infty}^{\alpha,\alpha+\frac{n}{p}-1}}.
\]
From the last two estimates, we get%

\begin{equation}
\left\Vert G(t)u_{0}\right\Vert _{X}\leq C\left\Vert u_{0}\right\Vert
_{\dot{B}W\dot{K}_{p,q,\infty}^{\alpha,\alpha+\frac{n}{p}-1}}.
\label{eq:Estimativa termo linear em X}%
\end{equation}

Take $0<\epsilon<1/4K$ and $0<\delta<\epsilon/C$ where $C$ is as in
(\ref{eq:Estimativa termo linear em X}). It follows from
\eqref{eq:Estimativa termo linear em X} and
\eqref{eq:Estimativa termo naolinear em X} that%
\begin{align*}
\left\Vert \Psi(u)\right\Vert _{X}  &  \leq\left\Vert G(t)u_{0}\right\Vert
_{X}+\left\Vert B(u,u)\right\Vert _{X}\\
&  \leq C\left\Vert u_{0}\right\Vert _{\dot{B}W\dot{K}_{p,q,\infty}%
^{\alpha,\alpha+\frac{n}{p}-1}}+K\left\Vert u\right\Vert _{X}\left\Vert
v\right\Vert _{X}\\
&  \leq2\epsilon
\end{align*}
So, $\Psi$ is well-defined, moreover for $u,v\in\bar{B}(0,2\epsilon)$ we have that%

\begin{align}
\left\Vert \Psi(u)-\Psi(v)\right\Vert _{X}  &  =\left\Vert
B(u-v,u)+B(v,u-v)\right\Vert _{X}\nonumber\\
&  \leq K\left\Vert u-v\right\Vert _{X}\left\Vert u\right\Vert _{X}%
+K\left\Vert v\right\Vert _{X}\left\Vert u-v\right\Vert _{X}\nonumber\\
&  \leq4K\epsilon\left\Vert u-v\right\Vert _{X}. \label{aux-cont-dependence-1}%
\end{align}
Since $4K\epsilon<1$, we get that $\Psi$ is a contraction and then this part
is concluded by the Banach fixed-point theorem. Notice that the continuous
dependence with respect to the initial data $u_{0}$ follows from estimates
(\ref{eq:Estimativa termo linear em X}) and (\ref{aux-cont-dependence-1}).

\

\textbf{Time-weak continuity at $t=0$.} The proof of the weak-$\ast$
convergence follows from the two following lemmas.

\

The first one is due to Kozono and Yamazaki \cite[pg. 989.]{KoYa}.

\begin{lemma}
For every real number $s$ and $u_{0}\in\dot{B}_{\infty,\infty}^{s},$ we have
$G(t)u_{0}\overset{\ast}{\rightharpoonup}u_{0}$ in $\dot{B}_{\infty,\infty
}^{s}$ as $t\rightarrow0^{+}.$
\end{lemma}

\fin

The second one is concerned with the weak-convergence of the bilinear term
$B(u,u)$ and it concludes the proof.

\begin{lemma}
Let $v\in X.$ We have that $B(v,v)(t)$ converges to $0$ in the weak-$\ast$
topology of $\dot{B}_{\infty,\infty}^{-1}$ as $t\rightarrow0^{+}$.
\end{lemma}

\textbf{Proof.} Let $\phi\in\dot{B}_{1,1}^{1}$ and $\epsilon>0$ an arbitrary
number. We can choose $\tilde{\phi}\in\mathcal{S}$ such that $\left\Vert
\phi-\tilde{\phi}\right\Vert _{\dot{B}_{1,1}^{1}}<\epsilon$. Then we have that%

\begin{align}
\left\vert \left\langle B(v,v)(t),\phi-\tilde{\phi}\right\rangle \right\vert
&  \leq\left\Vert B(v,v)(t)\right\Vert _{\dot{B}_{\infty,\infty}^{-1}%
}\left\Vert \phi-\tilde{\phi}\right\Vert _{\dot{B}_{1,1}^{1}}\nonumber\\
&  \leq C\left\Vert B(v,v)(t)\right\Vert _{\dot{B}W\dot{K}_{p,q,r}%
^{\alpha,\alpha+n/p-1}}\left\Vert \phi-\tilde{\phi}\right\Vert _{\dot{B}%
_{1,1}^{1}}\leq K\left\Vert v\right\Vert _{X}^{2}\epsilon\leq C\epsilon.
\label{Aux. Cfraca1}%
\end{align}
On the other hand,%

\begin{align}
\left\vert \left\langle B(v,v)(t),\tilde{\phi}\right\rangle \right\vert  &
\leq\int_{0}^{t}\left\vert \left\langle G(t-\tau)\mathbb{P}%
\mbox{{\rm div}}\left[  v\otimes v\right]  (\tau),\tilde{\phi}\right\rangle
\right\vert d\tau\leq\int_{0}^{t}\left\vert \left\langle \mathbb{P}%
\mbox{{\rm div}}\left[  v\otimes v\right]  (\tau),G(t-\tau)\tilde{\phi
}\right\rangle \right\vert d\tau\nonumber\\
&  \leq\int_{0}^{t}\left\Vert \mbox{{\rm div}}\left[  v\otimes v\right]
(\tau)\right\Vert _{\dot{B}_{\infty,\infty}^{-1-2\alpha-n/p}}\left\Vert
G(t-\tau)\tilde{\phi}\right\Vert _{\dot{B}_{1,1}^{1+2\alpha+n/p}}%
d\tau\nonumber\\
&  \leq C_{\tilde{\phi}}\int_{0}^{t}\left\Vert \left[  v\otimes v\right]
(\tau)\right\Vert _{\dot{B}_{\infty,\infty}^{-2\alpha-n/p}}d\tau\leq
C_{\tilde{\phi}}\int_{0}^{t}\left\Vert \left[  v\otimes v\right]
(\tau)\right\Vert _{W\dot{K}_{p,q}^{2\alpha}}d\tau\nonumber\\
&  \leq C_{\tilde{\phi}}\int_{0}^{t}\tau^{\left[  -\frac{1}{2}+\left(
\frac{\alpha}{2}+\frac{n}{4p}\right)  \right]  \cdot2}\tau^{\left[  \frac
{1}{2}-\left(  \frac{\alpha}{2}+\frac{n}{4p}\right)  \right]  \cdot
2}\left\Vert v(\tau)\right\Vert _{W\dot{K}_{2p,2q}^{\alpha}}^{2}%
d\tau\nonumber\\
&  \leq C_{\tilde{\phi}}\left\Vert v\right\Vert _{X}^{2}\int_{0}^{t}%
\tau^{-1+\alpha+\frac{n}{2p}}d\tau\leq C_{\tilde{\phi}}\left\Vert v\right\Vert
_{X}^{2}t^{\alpha+\frac{n}{2p}}. \label{Aux. Cfraca2}%
\end{align}
From \eqref{Aux. Cfraca1} and \eqref{Aux. Cfraca2}, we obtain
\begin{align*}
0\leq\mathop{\mbox{limsup}}\limits_{t\rightarrow{0^{+}}}\left\vert
\left\langle B(v,v)(t),\phi\right\rangle \right\vert  &  \leq
\mathop{\mbox{limsup}}\limits_{t\rightarrow{0^{+}}}\left\vert \left\langle
B(v,v)(t),\phi-\tilde{\phi}\right\rangle \right\vert
+\mathop{\mbox{limsup}}\limits_{t\rightarrow{0^{+}}}\left\vert \left\langle
B(v,v)(t),\tilde{\phi}\right\rangle \right\vert \\
&  \leq C\epsilon+0.
\end{align*}
Since $\epsilon>0$ is arbitrary, we conclude that
$\mathop{\lim}\limits_{t\rightarrow{0^{+}}}\left\vert \left\langle
B(v,v)(t),\phi\right\rangle \right\vert =0.$ Now, using that $\phi\in\dot
{B}_{1,1}^{1}$ is arbitrary, we get the desired convergence. \fin

\end{document}